\documentclass[reqno,sort&compress]{siamonline220329}

\usepackage{graphicx} 

\usepackage[sort]{cite} 
\usepackage{geometry}
\usepackage{comment}
\usepackage[english]{babel}

\usepackage[utf8]{inputenc}

\usepackage{hhline}
\usepackage{bbm}
\usepackage{multirow}
\usepackage{tabularx}
\usepackage{amsmath,amssymb}
\usepackage{siunitx}
\usepackage{placeins}
\usepackage{longtable}
\usepackage{graphicx}
\usepackage{tikz}

\usepackage{color}

\usepackage{algorithm}
\usepackage[noend]{algpseudocode}

\newcommand*\Let[2]{\State #1 $\gets$ #2}

\mathchardef\ordinarycolon\mathcode`\:
\mathcode`\:=\string"8000
\begingroup \catcode`\:=\active
\gdef:{\mathrel{\mathop\ordinarycolon}}
\endgroup
\usepackage{tikz}
\usetikzlibrary{calc,arrows}

\newcommand{\norm}[2]{\left\| #1 \right\|_{#2}}

\newcommand{\pdiff}[2]{\dfrac{\partial{#1}}{\partial{#2}}}

\newcommand{\diff}[2]{\dfrac{d{#1}}{d{#2}}}

\newcommand{\bv}[1]{{\bf #1}}
\newcommand{\bs}[1]{\boldsymbol{#1}}
\newcommand{\by}{{\mathbf{y}}}

\newcommand{\IR}{{\mathbb{R}}}

\renewcommand{\leq}{\leqslant}
\renewcommand{\geq}{\geqslant}

\newcommand{\Kxi}{K_{x}^{i}}
\newcommand{\Kyj}{K_{y}^{j}}

\newcommand{\nss}[1]{n_{#1}}
\newcommand{\rss}[1]{r_{#1}}

\newcommand{\half}{\frac12}

\newcommand{\wLi}[1]{w^{L,i}_{#1}}
\newcommand{\wRi}[1]{w^{R,i}_{#1}}
\newcommand{\RANK}{\textrm{rank}}

\newcommand{\mat}[1]{\mathcal{#1}}
\newcommand{\ten}[1]{\mathcal{#1}}
\newcommand{\tenb}[1]{\bs{\mathcal{#1}}}

\newcommand{\idx}[1]{#1}

\title{The Tensor-Train Stochastic Finite Volume Method for
  Uncertainty Quantification}

\author{
  Steven Walton\thanks{ T-5 Group, Theoretical Division, Los Alamos National Laboratory \email{[swalton,tokareva, gmanzini]@lanl.gov}}
  \and Svetlana Tokareva\footnotemark[1]
  \and Gianmarco Manzini\footnotemark[1]
}

\date{\today}

\begin{document}

\maketitle

\begin{abstract}
  The stochastic finite volume method offers an efficient one-pass
  approach for assessing uncertainty in hyperbolic conservation
  laws.
  Still, it struggles with the curse of dimensionality when dealing
  with multiple stochastic variables.
  We introduce the stochastic finite volume method within the
  tensor-train framework to counteract this limitation.
  This integration, however, comes with its own set of difficulties,
  mainly due to the propensity for shock formation in hyperbolic
  systems.
  To overcome these issues, we have developed a tensor-train-adapted
  stochastic finite volume method that employs a global WENO
  reconstruction, making it suitable for such complex systems.
  This approach represents the first step in designing tensor-train
  techniques for hyperbolic systems and conservation laws involving
  shocks.
\end{abstract}

\begin{keywords}
  low-rank approximation,
  tensor-train decomposition,
  stochastic finite volume method,
  uncertainty quantification,
  hyperbolic systems,
  conservation laws,
  discontinuous solutions,
  shocks
\end{keywords}

\begin{AMS}
  15A69, 65C20, 76M12, 35L65, 65F10, 65F50
\end{AMS}


\raggedbottom

\section{Introduction}

The field of uncertainty quantification has exploded in recent years,
cf.~\cite{DElia-Gunzburger-Rozza:2020}.
Parametrized and stochastic partial differential equations are
particularly interesting in this burgeoning field.
Numerical methods for uncertainty quantification of solutions to
Partial Differential Equations (PDEs) often fall into two camps:
intrusive and non-intrusive, cf.~\cite{zbMATH06580444}.
A recent method, the Stochastic Finite Volume (SFV)
method~\cite{SFVSpringer2014,AbgrallTokareva-2017,BAR1}, has emerged
as what we are calling a weakly intrusive method as it requires very
little modification of the existing code and does not alter properties
of the system under consideration, e.g., hyperbolic systems remain
hyperbolic.
In contrast to Polynomial Chaos Expansion
(PCE)~\cite{Knio-LeMaitre,GrotePainTeckentrupVandenHof2015}, which
employs global polynomials, the SFV method represents the stochastic
solution as a piecewise constant function across the stochastic domain
discretization.
Higher-order approximations can be obtained in every mesh cell through
local piecewise polynomial reconstruction techniques such as the WENO
method.
The SFV method does not require continuity between the discrete
elements, thus providing a compact discretization stencil and
allowing for efficient parallelization.
While the SFV method has notable advantages over Monte Carlo
simulations in terms of accuracy, and PCE methods in its ability to
preserve hyperbolicity and monotonicity, it faces computational
challenges.
Specifically, when tackling problems with a high number of stochastic
parameters, the SFV method quickly becomes computationally prohibitive
due to the \emph{curse of dimensionality}~\cite{Bellman:1961}.

A particularly effective strategy to overcome the curse of
dimensionality, an exponential increase in complexity with the number
of dimensions of the problem, consists in combining tensor network
formats with low-rank approximations, see, e.g.,
\cite{Grasedyck-2013}.
In recent years, tensor decompositions have experienced a surge in
applications across various
fields~\cite{Kolda2009TensorDecompositions}.
Of the many tensor decomposition techniques available, the
tensor-train (TT) format~\cite{OSELEDETS201070,Oseledets-TT-2011}
stands out due to its robustness and beneficial properties.
Precisely, it exhibits a linear relationship between complexity and
dimensionality and offers provable tensor rank and accuracy bounds,
thus effectively separating the dimensions of the problem.
A TT-format reformulation of the SFV method makes it possible to
overcome the curse of dimensionality when a low-rank approximate
solution to the stochastic equations exists.

The transformation of numerical discretizations into the TT framework
has a track record of success in solving ordinary and partial
differential equations, see,
e.g.,\cite{MANZINI2023615,Truong-etal:2023,Kazeev2012LowRankEQ,Kazeev2016QTTfiniteelementAF,Kazeev2018QuantizedTF,Khoromskij2009QuanticsTTAO,Khoromskij2011QTTAO,Khoromskij2012NumericalAO,Dolgov-Khoromskij:2014,Dolgov:2014,Dolgov2012FastSO}
with numerous achievements related to uncertainty quantification,
cf.~\cite{Dolgov_PCE_2015,Zhang_TT_ANOVA_2015,gruhlke2022lowrank,wang2023tensor}.
Despite the extensive body of work in this area, there appears to be a
gap in the literature: none of the existing studies have addressed the
numerical treatment of stochastic solutions with space discontinuities
like shocks, which are an inherently important characteristic of
hyperbolic problems.
Our work aims to bridge this gap by developing an SFV method within
the TT framework that is specifically designed to address these
complex features.
To this end, we reformulate the cell-wise weighted essentially
non-oscillatory (WENO) finite volume schemes~\cite{Shu} within the
tensor-train framework.
We achieve this task by creating global reconstruction matrix
operators that facilitate the direct application of polynomial
reconstruction to the TT cores, allowing the cell-averaged states to
be effectively processed.
Additionally, we introduce global quadrature reconstruction matrices
that similarly apply directly to the decomposition cores.
These reconstruction matrices are instrumental in accurately rendering
discontinuous solution profiles within the TT format while avoiding
the generation of spurious oscillations at discontinuities.
The global reconstruction matrices are {\emph {applicable
  beyond uncertainty quantification problems to general hyperbolic
  systems}}.
This topic will be the subject of forthcoming studies.

\subsection*{Outline}
The structure of this article is as follows.
Section \ref{sec::SFV} offers a concise overview of the SFV method, in which we reformulate the
one-dimensional WENO reconstruction algorithm as a global matrix
operator.
In Section \ref{sec::TT-SFV-method}, we discuss a few fundamental
components of the tensor-train format and we present a comprehensive
exposition of the SFV methodology within the TT framework.
Section \ref{sec::tests} showcases a series of numerical experiments
demonstrating the method's accuracy, stability, and efficiency,
particularly as the number of stochastic dimensions grows.
This section includes performance evaluations of the TT-SFV algorithm
on a Sod shock tube problem under varying stochastic dimensions,
alongside a comparative analysis of our Python implementation against
a traditional C++ SFV code.
Section \ref{sec::conclusion} finally discusses the current
limitations of our TT-SFV approach and proposes potential topics for
future research.



\section{Stochastic Finite Volume Method}
\label{sec::SFV}

In this section, we provide a comprehensive explanation of the SFV
method for hyperbolic PDE systems.  We restrict the present discussion to one physical dimension.  The extension of ENO/WENO finite volume schemes to more than one physical dimension can be found in \cite{Shu}.  Similarly, the extension to more than one physical dimension of the SFV method is identical, provided the stochastic variables are taken into account as discussed below. In section \ref{sec::TT-SFV-method}, describing the reformulation of the SFV method in TT-format, we will no longer restrict ourselves to a 1-dimensional physical space. 

Let $p$ be a positive integer, $t\in(0, T]=I$ with $T>0$, $x\in\big[x_L, x_R\big]=D_x\subset \IR^n$ for $n=1$ and $\omega\in\Omega$,
which is the set of all possible outcomes of the model parameter's
uncertainties. We consider the hyperbolic system of first-order PDEs in one space
dimension
\begin{subequations}
  \begin{alignat}{2}
    \pdiff{\bv{U}}{t}+\pdiff{\bv{F}(\bv{U})}{x} & = \bv{S}(\bv{U},\omega),  &&\qquad (t, x, \omega)\in I\times D_x \times \Omega, \label{eq:1D-PDE} \\[0.125em]
    \bv{U}(0,x,\omega)                          & = \bv{U}_0  (x,\omega),   &&\qquad (t, x, \omega)\in \{0\}\times D_x \times \Omega,       \label{eq:1D-IC}  \\[0.125em]
    \bv{U}(t,x_L,\omega)                        & = \bv{U}_B^L(t,\omega),    &&\qquad (t, x, \omega)\in I\times \{x_L\} \times \Omega,                 \label{eq:1D-BC-L}\\[0.125em]
    \bv{U}(t,x_R,\omega)                        & = \bv{U}_B^R(t,\omega).    &&\qquad (t, x, \omega)\in I\times \{x_R\} \times \Omega,                 \label{eq:1D-BC-R}
  \end{alignat}
\end{subequations}
where $\bv{U}=(U_1,U_2,\ldots,U_p)^T$ is the vector of $p$ unknown
conserved variables $U_{\ell}(t,x,\omega)$, $\ell=1,2,\ldots,p$, and
$\bv{F}=(F_1,F_2,\ldots,F_p)^T$ is the corresponding flux vector.
Each $U_{\ell}(x,t,\omega)$ is a function of time $t$, space $x$ and a
set of stochastic independent variables $\omega$.
The independent variable $\omega$ introduces uncertainty through the
initial and boundary conditions, ie., \eqref{eq:1D-IC} and
\eqref{eq:1D-BC-L}-\eqref{eq:1D-BC-R} respectively, and the forcing
term in the right-hand side of Eq. \eqref{eq:1D-PDE}.
Each $F_{\ell}(\bv{U})$, the $\ell$-th component of the flux vector
$\bv{F}$, depends on the variable vector $\bv{U}(t,x,\omega)$, thus
making the system of equations~\eqref{eq:1D-PDE} nonlinear.
This system is hyperbolic if the $p\times p$-sized, Jacobian matrix
$\partial\bv{F}\slash{\partial\bv{U}}=\big(\partial F_i\slash{\partial
  U_j}\big)_{i,j=1,\ldots,p}$ has $p$ real eigenvalues and $p$
linearly independent eigenvectors for every $\bv{U}$ in the domain of
interest, cf.~\cite{Godlewski-Raviart,LeVeque,evans2010partial}.

Let $m$ be a positive integer, and assume we can parametrize $\omega$
using the random variable $\by=Y(\omega)\in D_{\by} \subset\IR^m$
with $\mathbb{E}[\by]  = \int_{D_{\by}} \by \mathrm d \mu(\by) < \infty$, where $\mu(\by)$ is the probability density function on $D_\by$. We assume that $\mu(\by)$ is bounded on $D_\by$ or that the error induced by truncation onto $D_\by$ is small in the case $\mu(\by)$ is unbounded.  See \cite{sfv-adaptive-SIAM} and \cite{HertyKolbMuller2023} for further relevant discussion.
Then the initial condition, the boundary data and the source term take the
following form: 
\begin{subequations}
  \begin{align}
    \bv{U}_0(x,\by)    &:= \bv{U}(0,x,   Y(\omega)), \label{eq:1D-IC-par} \\[0.25em]
    \bv{U}_B^L(t,\by)  &:= \bv{U}(t,x_L, Y(\omega)), \label{eq:1D-BC0-par} \\[0.25em]
    \bv{U}_B^R(t,\by)  &:= \bv{U}(t,x_R, Y(\omega)), \label{eq:1D-BC1-par} \\[0.25em]
    \bv{S}(\bv{U},\by) &:= \bv{S}(\bv{U},Y(\omega)),  \label{eq:1D-SRC-par}
  \end{align}
\end{subequations}
where, with some minor abuse of notation, we still use the symbols ``$\bv{U}(\cdot,\cdot,\cdot)$'' for the
conserved quantities.
After this transformation, the parametric form of the stochastic PDE
system \eqref{eq:1D-PDE}-\eqref{eq:1D-BC-R} becomes
\begin{subequations}
  \vspace{-0.8\baselineskip}
  \begin{alignat}{2}
    \pdiff{\bv{U}}{t}+\pdiff{\bv{F}}{x} &= \bv{S}(\bv{U},\by), &&\qquad (t, x,\by) \in I \times D_x \times D_{\by}, \label{eq:1D-PDE-par}  \\[0.125em]
    \bv{U}(0,x,  \by)                   &= \bv{U}_0(x,\by),    &&\qquad (t, x,\by) \in \{0\} \times D_x \times D_{\by},         \label{eq:1D-IC-par-1} \\[0.125em]
    \bv{U}(t,x_L,\by)                   &= \bv{U}_B^L(t,\by),   &&\qquad (t, x,\by) \in I \times \{x_L\} \times D_{\by},                  \label{eq:1D-BC0}      \\[0.125em]
    \bv{U}(t,x_R,\by)                   &= \bv{U}_B^R(t,\by),    &&\qquad (t, x,\by) \in I \times \{x_R\} \times D_{\by}.                  \label{eq:1D-BCL}
  \end{alignat}
\end{subequations}
We can now consider \eqref{eq:1D-PDE-par}--\eqref{eq:1D-BCL} as a
system of partial differential equations for the time-dependent state
vector $\bv{U}(t,x,\bv{y})$ as a function of the $d$-dimensional
independent variables $(x,\by)\in D_x\times
D_{\bv{y}}\subset\mathbb{R}^n\times\mathbb{R}^{m}$ with $d=n+m$ and discretize it
with the multidimensional finite volume method.

The physical domain is partitioned into a collection of non-overlapping cells, $K_x^i\subset D_x$, with $\bigcup_{i}K_x^i = D_x$. Likewise, we partition the stochastic
domain into a finite set of closed, equi-spaced, $m$-dimensional
intervals $K_\by^j\subset D_\by$, where $j$ is the $m$-sized multi-index
$j=(j_{\alpha})_{\alpha=1,\ldots,m}$.
The Cartesian product of the partitions of the physical and stochastic
domains forms the multidimensional mesh, i.e.,
$\{\Kxi\}\times\{K_\by^j\}$.

Then, we introduce the cell-average of the state vector
$\bv{U}(t,x,\by)$ over the $ij$-th $d$-dimensional cell $\Kxi\times K_\by^j$
as
\begin{align}
  |K_\by^j|\bar{\bv{U}}_{ij}(t) = \mathbb{E}_{\by_j}\Bigg[\frac{1}{\Delta x_i}\int_{\Kxi}\bv{U}(t,x,\by)\,dx\Bigg].
  \label{eq:Ubar:def}
\end{align}
$\Delta x_i = |\Kxi|=\int_{\Kxi}\mathrm{d}x$ is the length of $\Kxi$
in the 1-dimensional physical space. We will use $|\Kxi|$ rather than $\Delta x$ when $ n > 1$ in section \ref{sec::TT-SFV-method}.
The expectation $\mathbb{E}_{y_j}[\,\cdot\,]$ is taken over the
cell $K_\by^j$ with respect to the measure
$\mu(\by)$, and 
$|K_\by^j|=\int_{K_\by^j}\mathrm{d}\mu(\by)$ is the measure of the cell $K_\by^j$.
Integrating over the control volume and applying the divergence theorem, we obtain (the still exact) formulation
\begin{equation}\label{eqn::sfv_exact}
  \Delta x_i|K_\by^j|\diff{\bv{\bar{U}}_{ij}}{t}
  + \mathbb{E}_{\by_j}\Big[\bv{F}\big(\bv{U}(t,x_{i+\half},\by)\big) - \bv{F}\big(\bv{U}(t,x_{i-\half},\by)\big) \Big]
  = \mathbb{E}_{\by_j}\Bigg[\frac{1}{\Delta x_i}\int_{\Kxi}\bv{S}\big(\bv{U}(t,x,y)\big)\,dx\Bigg].
\end{equation}
To compute the flux integrals, we apply a suitable quadrature rule at the interfaces $x_{i\pm\frac12}$ 
given by the collection of nodes and weights $(x_{i\pm\frac12}, \by_{\bs{m}}, \omega_{\bs{m}})$ with $\bs{m}$ a multi-index labeling the quadrature points on the interfaces of the
stochastic cell $\Kyj$. In accordance with the formulation of the Godunov discretization
method~\cite{LEV1}, we introduce the numerical flux $\hat{\bv{F}}\big(
  \bv{U}^{\bs m,L}_{i\pm\half,j},\bv{U}^{\bs m,R}_{i\pm\half,j}\big)$ which depends on $\bv{U}^{\bs m,L}_{i\pm\half,j}$ and $\bv{U}^{\bs m,R}_{i\pm\half,j}$ where we use the standard convention of indicating the ``left'' and ``right'' interface values at $x_{i\pm\frac12}$ with the letters $L$ and $R$.

The numerical flux $\hat{\bv{F}}(\cdot,\cdot)$ can be any flux formula
from the shock-capturing methodology, e.g., Godunov, Lax-Friedrichs,
HLLC flux, etc.
A specific example of a numerical flux, and the one we use for every numerical test in section \ref{sec::tests}, is the Lax-Friedrichs
flux,
\begin{equation}
  \hat{\bv{F}}\big(\bv{U}^{\bs m,L}_{i\pm\half,j},\bv{U}^{\bs m,R}_{i\pm\half,j}\big)
  = \dfrac{1}{2}\Big(\bv{F}(\bv{U}^{\bs m,L}_{i\pm\half,j})
  + \bv{F}(\bv{U}^{\bs m,R}_{i\pm\half,j})\Big)
  - \dfrac{\nu}{2}\Big(\bv{U}^{\bs m,R}_{i\pm\half,j}-\bv{U}^{\bs m,L}_{i\pm\half,j}\Big),
  \label{eq:laxfriedrichsflux}
\end{equation}
where $\nu$ is a user-defined, numerical viscosity coefficient, and
$\bv{U}_{L}$ and $\bv{U}_{R}$ denote the left and right interface
vales.  We take $$\nu = \max\Bigg\{ \Big|\mathbf F^\prime \big(\bv{U}^{\bs m,L}_{i\pm\half,j}\big) |, |\mathbf F^\prime \big(\bv{U}^{\bs m,R}_{i\pm\half,j}\big) \Big|\Bigg\},$$ i.e. the Rusanov or (approximate) local Lax-Friedrich's numerical flux \cite{KuzminHajduk.2023}.

Evaluation of the numerical flux at the quadrature points leads to the following definition
\begin{align}
  \bar{\bv{F}}_{i\pm\half,j}(t) := \sum_{\bs m}
  \hat{\bv{F}}\big(
  \bv{U}^{\bs m,L}_{i\pm\half,j},\bv{U}^{\bs m,R}_{i\pm\half,j}\big)\,\mu(\by_j^{\bs m})\omega_{\bs m}
  \approx \mathbb{E}_{\by_j}\Big[ \bv{F}\big(\bv{U}(t, x_{i\pm\half}, \by_j)\big) \Big].
  \label{eq:numflux:def}
\end{align}
Similarly, we write the approximation of the source term as
\begin{align}
  \bar{\bv{S}}_{ij}(t) := \sum_{\bs q}
  \bv{S}\big(
  \bv{U}(t,x_{\bs{q}},\by_j^{\bs{q}})\big)\mu(\by_j^{\bs q})\omega_{\bs q}\approx \mathbb{E}_{\by_j}\Bigg[\frac{1}{\Delta x_i}\int_{\Kxi}\bv{S}\big(\bv{U}(t,x,\by)\big)\,\mathrm d x\Bigg],
  \label{eq:Sbar:def}
\end{align}
where the multi-index $\bs q$ labeling the quadrature nodes and weights in $\Kxi\times K_\by^j$.  An identical quadrature is applied to \eqref{eq:Ubar:def} and, abusing notation, we denote this approximation by $\bv{\bar{U}}_{ij}$ also.  This equation and both definitions~\eqref{eq:Ubar:def} and
\eqref{eq:numflux:def} must be intended componentwise, i.e., they
independently hold for each component $U_{\ell}$ of $\bv{U}$ and
$F_{\ell}$ of $\bf{F}$, $\ell=1,2,\ldots,p$.

Then, the semi-discrete stochastic finite volume method reads simply as
\begin{equation}
  \frac{\mathrm{d}\bar{\bv{U}}_{ij}(t)}{\mathrm d t}
  + 
  \frac{1}{\Delta x}
  \bigg[
    \bar{\bv{F}}_{i+\half,j}(t) - \bar{\bv{F}}_{i-\half,j}(t)
    \bigg] = \bar{\bv{S}}_{ij}(t),
    \label{eqn::general_SFV}
\end{equation}
with the initial condition $\bv{U}_{ij}(0)$ given by 
\begin{equation}\label{eqn::init_Ubar}
  \bv{\bar{U}}_{ij}(0) = \sum_{\bs q}
  \bv{{U}}(0,x_{\bs{q}},\by_j^{\bs{q}})\mu(\by_j^{\bs q})\omega_{\bs q} \approx \mathbb{E}_{\by_j}\Bigg[\frac{1}{\Delta x_i}\int_{\Kxi}\bv{{U}}(0,x,\by)\big)\,\mathrm d x\Bigg].
\end{equation}

\medskip
Let $I=[0,T]$ be the time integration interval from $t=0$ up to
the final time $T>0$, and consider a uniform partition into $K$
sub-intervals $\big[t^{k},t^{k+1}\big]$, $k=0,\ldots,K$, with
$0=t^{0}<\ldots<t^{k}<t^{k+1}<\ldots<t^{K}=T$.
We obtain a first-order, explicit discretization in time using the
forward Euler time-marching scheme, which reads
\begin{equation}
  \bar{\bv{U}}_{ij}^{k+1}
  = \bar{\bv{U}}_{ij}^{k}
  - \frac{\Delta t}{\Delta x} \bigg[\bar{\bv{F}}^{k}_{i+\frac12,j} - \bar{\bv{F}}^{k}_{i-\frac12,j}\bigg]
  + \Delta t\bar{\bv{S}}_{ij}^k
  \label{eqn::general:update}
\end{equation}
with the numerical flux terms
$\bar{\bv{F}}^{k}_{i\pm\frac12,j}:=\bar{\bv{F}}_{i+\half,j}(t^k)$
from~\eqref{eq:numflux:def} and the source term
$\bar{\bv{S}}_{ij}^k=\bar{\bv{S}}_{ij}(t^k)$ evaluated using
$\bar{\bv{U}}_{ij}^{k} := \bar{\bv{U}}_{ij}(t^k)$, and with the initial solution
$\bv{U}_{ij}^{0}$ given again by \eqref{eqn::init_Ubar}.

It is well known that higher-order accurate, strong stability preserving (SSP) explicit time discretizations can be obtained by
using~\eqref{eqn::general:update} as a building block (see \cite{Gottlieb2009, Shu-1987} among others).  In practice, we use either a second- or third-order SSP explicit scheme.  Which time integration scheme is used is indicated in each of the numerical tests.

To elaborate on the interface reconstruction and ease notation, we further simplify the above discussion and assume that $n=m=1$ or $d=2$ for the remainder of the subsection. Hence, we obtain a two-dimensional, time-dependent problem in $x$ and $y$.
We again take a uniform grid in both
physical and stochastic variables $x$ and $y$, with nodes
$x_{i\pm\half}$ and $y_{j\pm\half}$, respectively spaced by $\Delta x$ as before and
now the measure of cell $(i,j)$ in the $y$-direction
defined by 
\begin{align*}
  |\Delta y| = \int\limits_{y_{j-\half}}^{y_{j+\half}} \mu(y)\,\mathrm dy.
\end{align*}
The $(i,j)$-cell delimited by this partition is the rectangle
$\big[x_{i-\half},x_{i+\half}\big]\times\big[y_{j-\half},y_{j+\half}\big]$.
We arrive at simplified definitions of the cell average, numerical flux, source term, and initial condition.  The cell average of the unknown vector $\bv{U}(t,x,y)$ over the
$(i,j)$-th cell is defined as
\begin{align*}
  \bv{\bar{U}}_{ij}(t)
  = \dfrac{1}{\Delta x |\Delta y|}
  \int\limits_{x_{i-\half}}^{x_{i+\half}}\int\limits_{y_{j-\half}}^{y_{j+\half}} \bv{U}(t,x,y)\mathrm d\mu(y)\mathrm dx.
\end{align*}
The analogue of equation \eqref{eqn::sfv_exact} is
\begin{multline}
  \Delta x |\Delta y| \diff{\bv{\bar{U}}_{ij}}{t}
  + \int\limits_{y_{j-\half}}^{y_{j+\half}} \Big( \bv{F}\big(\bv{U}(t,x_{i+\half},y)\big)\mu(y) - \bv{F}\big(\bv{U}(t,x_{i-\half},y)\big)\Big)\mathrm d \mu(y)
  = \int\limits_{x_{i-\half}}^{x_{i+\half}}\int\limits_{y_{j-\half}}^{y_{j+\half}} \bv{S}(\bv{U},y)\mathrm d\mu(y)\,dx,
\end{multline}
where we now have a simpler 1-dimensional quadrature rule
at the interface $x_{i+\half}$ with weights and nodes
$\big\{(w_m,x_{i+\half},y_j^m)_{m=1,\ldots,M}\big\}$.

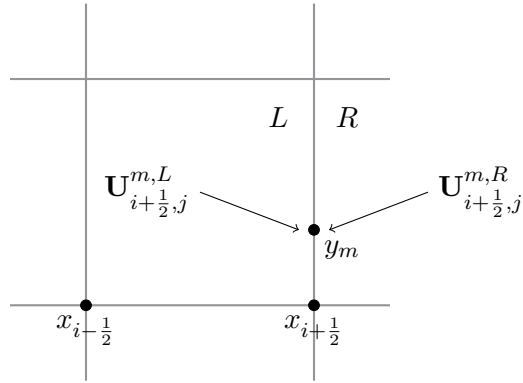
\begin{figure}[h!]
  \begin{center}
    \begin{tikzpicture}
      \draw[gray, thick] (0,0) -- (5,0);
      \draw[gray, thick] (0,3) -- (5,3);
      \draw[gray, thick] (1,-1) -- (1,4);
      \draw[gray, thick] (4,-1) -- (4,4);
      \filldraw[black] (1,0) circle (2pt) node[anchor=north] {$x_{i-\half}$};
      \filldraw[black] (4,0) circle (2pt) node[anchor=north] {$x_{i+\half}$};
      \filldraw[black] (4,1) circle (2pt) node[anchor=north west] {$y_m$};
      \draw[->] (2.5,1.5) -- (3.8,1);
      \draw[->] (5.5,1.5) -- (4.2,1);
      \draw (2.5,1.5) node[anchor=east] {${\bv{U}}^{m,L}_{i+\half,j}$};
      \draw (5.5,1.5) node[anchor=west] {${\bv{U}}^{m,R}_{i+\half,j}$};
      \draw (3.8,2.5) node[anchor=east] {$L$};
      \draw (4.15,2.5) node[anchor=west] {$R$};
    \end{tikzpicture}
  \end{center}
  \caption{Illustration of reconstruction points for the SFV method;
    ``L'' and ``R'' label the ``left'' and ``right'' sides of the cell
    interface located at $x_{i+\half}$, and
    ${\bv{U}}^{m,L}_{i+\half,j}$ and
    ${\bv{U}}^{m,R}_{i+\half,j}$ are the reconstructed values at
    the quadrature node with coordinates $(x_{i+\half},y_m)$. }
  \label{Fig:SFV-weakly-intrusive}
\end{figure}

We again use the definition \eqref{eq:laxfriedrichsflux} for the numerical flux
$\hat{\bv{F}}\big({\bv{U}}_{i+\half,j}^{m,L}(t),{\bv{U}}_{i+\half,j}^{m,R}(t)\big)$,
which depends on ${\bv{U}}_{i+\half,j}^{m,L}(t)$, and ${\bv{U}}_{i+\half,j}^{m,R}(t)$ as before. See Fig.~\ref{Fig:SFV-weakly-intrusive}.
We obtain these values by taking the trace at the interface $j$ of a
reconstructed approximate solution, a piecewise polynomial (see the next subsection), inside the
two adjacent cells $(i,j)$ and $(i+1,j)$.
Then, evaluating these polynomials at the quadrature nodes
$(x_{i+\half},y_j^m)$ on the interface $j$ yields
\begin{align}
  \int\limits_{y_{j-\half}}^{y_{j+\half}} \bv{F}\big(\bv{U}(t,x_{i+\half},y)\big)\mu(y)\,dy
  \approx \sum\limits_{m=1}^M \hat{\bv{F}}\big({\bv{U}}_{i+\half,j}^{m,L}(t),{\bv{U}}_{i+\half,j}^{m,R}(t)\big)\,\mu(y_j^m)w_m
  =: \bar{\bv{F}}_{i+\half,j}(t).
  \label{eq:quadrature}
\end{align}
The remaining definitions for the source term and initial condition are identical, mutatis mutandis.  

In the next subsection, we describe the reconstruction procedure alluded to in the foregoing discussion. We will again assume that $n=1$ and describe the WENO reconstruction process on a scalar quantity $u$ which illustrates the process for one of our conserved quantities $U_\ell$  of $\bv{U}$, $\ell=1,2,\ldots,p$.  

\subsection{Global matrix formulation of WENO polynomial reconstruction}
\label{subsec::reconstruction}

To reformulate the SFV method in the tensor train format (which will
be the topic of the next section), we need a \emph{global matrix
representation of the polynomial reconstruction operator} that
provides the interface values that are used in the numerical flux
formula~\eqref{eq:quadrature}.
This reformulation is an alternative to the usual reconstruction
scheme that works element by element.

To ease the exposition, we consider the univariate reconstruction
provided by the WENO3 scheme, the generalization to other higher-order
WENO-type reconstructions being deemed as straightforward.
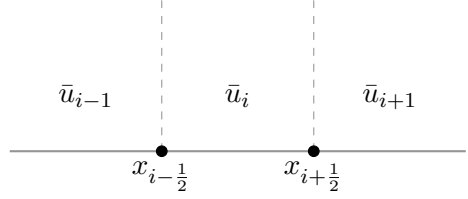
\begin{figure}[h!]
  \begin{center}
    \begin{tikzpicture}
      \draw[gray, thick] (0,0) -- (6,0);
      \draw[gray, dashed] (2,0) -- (2,2);
      \draw[gray, dashed] (4,0) -- (4,2);
      \filldraw[black] (2,0) circle (2pt) node[anchor=north] {$x_{i-\half}$};
      \filldraw[black] (4,0) circle (2pt) node[anchor=north] {$x_{i+\half}$};
      \draw (1,1) node[anchor=north] {${\bar{u}}_{i-1}$};
      \draw (3,1) node[anchor=north] {${\bar{u}}_i$};
      \draw (5,1) node[anchor=north] {${\bar{u}}_{i+1}$};
    \end{tikzpicture}
  \end{center}
  \caption{Reconstruction stencil of a scalar quantity for 1D WENO3.}
  \label{Fig:WENO3-1D-stencil}
\end{figure}

The reconstruction stencil of the WENO3 scheme consists of the cell
$i$ and its closest neighbors and is shown in
Fig.~\ref{Fig:WENO3-1D-stencil}; hence, we consider the stencil
centered at cell $i$ and also include the two cells indexed by
$i\pm1$.
We let $\bar{u}_i$ and $\bar{u}_{i\pm1}$ denote the cell averages of a
scalar, univariate quantity $u(x)$ over these cells.
Inside cell $i$, we consider the two linear polynomials
\begin{align*}
  p_0(x) = \bar{u}_i + \dfrac{\bar{u}_{i+1} - \bar{u}_i}{\Delta x} (x - x_i),\quad
  p_1(x) = \bar{u}_i + \dfrac{\bar{u}_i - \bar{u}_{i-1}}{\Delta x} (x - x_i).
\end{align*}
The reconstructed left and right value at the cell interface
$x_{i+\half}$ and $x_{i-\half}$ are
\begin{align*}
  u^L_{i+\half} = \wLi{0} p_0(x_{i+\half}) + \wLi{1} p_1(x_{i+\half}),\quad
  u^R_{i-\half} = \wLi{0} p_0(x_{i-\half}) + \wLi{1} p_1(x_{i-\half}), 
\end{align*}
where $\wLi{0}$, $\wLi{1}$, and $\wRi{0}$, $\wRi{1}$ are the
\emph{WENO interpolation weights} associated to cell $i$.
We emphasize that these weights can vary from cell to cell depending
on the smoothness of the function $u(x)$ that we reconstruct.
The smoothness of $u$ is detected by suitable smoothness indicators which approximate the variation of $u$ in the stencil,
cf.~\cite{Shu}.

Using the polynomial evaluations $p_0(x_{i\pm\half})$ and
$p_1(x_{i\pm\half})$ in the formulas for $u^L_{i-\half}$ and
$u^R_{i+\half}$ yields the formulas:
\begin{align*}
  u^R_{i-\half}
  &= \wRi{0} \bigg(\frac32\bar{u}_i - \frac12\bar{u}_{i+1}\bigg) + \wRi{1} \bigg(\frac12\bar{u}_{i-1}+\frac12\bar{u}_i\bigg) \\[0.5em]
  &= \frac12\wRi{1} \bar{u}_{i-1} + \bigg(\frac32\wRi{0} + \frac12\wRi{1}\bigg)\bar{u}_i - \frac12\wRi{0}\bar{u}_{i+1},\\[1em]
  u^L_{i+\half}
  &= \wLi{0} \bigg(\frac12\bar{u}_i + \frac12\bar{u}_{i+1}\bigg) + \wLi{1} \bigg(\frac32\bar{u}_{i}-\frac12\bar{u}_{i-1}\bigg) \\
  &= -\frac12\wLi{1} \bar{u}_{i-1} + \bigg(\frac12\wLi{0} + \frac32\wLi{1}\bigg)\bar{u}_i + \frac12\wLi{0}\bar{u}_{i+1}.
\end{align*}

We introduce the vectors
\begin{align}
  \hat{U} = (\bar{u}_1, \bar{u}_2, \ldots, \bar{u}_N)^T,\
  \hat{U}_L = \big(u_{\half}^L, u_{\frac32}^L, \dots, u_{N-\half}^L\big)^T, \
  \hat{U}_R = \big(u_{3/2}^R, \ldots, u_{N-\half}^R, u_{N+\half}^R\big)^T, 
\end{align}
and the matrices
\newcommand{\TABROW}[6]{ #1 & #2 & #3 & #4 & #5 & #6\\[0.5em]}
\begin{small}
  \begin{align}
    \mathbf{L} = 
    \begin{pmatrix}
      \TABROW{ \frac32w_0^{L,1} + \frac12w_1^{L,1}}{ -\frac12w_0^{L,1} }{ 0 }{ 0 }{ \ldots }{ 0 }
      \TABROW{ \frac12w_1^{L,2} }{ \frac32w_0^{L,2} + \frac12w_1^{L,2} }{ -\frac12w_0^{L,2} }{ 0 }{ \ldots }{ 0 }
      \TABROW{\vdots}{\vdots}{\vdots}{\vdots}{\vdots}{\vdots}
      \TABROW{ 0 }{ 0 }{ \ldots }{\frac12w_1^{L,N-1}}{\frac32w_0^{L,N-1} + \frac12w_1^{L,N-1}}{-\frac12w_0^{L,N-1}}
      \TABROW{ 0 }{ 0 }{ \ldots }{ 0 }{\frac12w_1^{L,N}}{\frac32w_0^{L,N} + \frac12w_1^{L,N}}
    \end{pmatrix}
    \nonumber
    \label{eqn::L}
  \end{align}
\end{small}
and
\begin{small}
  \begin{align}
    \mathbf{R} \!=\!\! 
    \begin{pmatrix}
      \TABROW{\frac12w_0^{R,1} + \frac32w_1^{R,1} }{ \frac12w_0^{R,1} }{ 0 }{ 0 }{ \ldots }{ 0 }
      \TABROW{-\frac12w_1^{R,2} }{ \frac12w_0^{R,2} + \frac32w_1^{R,2} }{ \frac12w_0^{R,2} }{ 0 }{ \ldots }{ 0 }
      \TABROW{ \vdots }{ \vdots }{  \vdots }{ \vdots }{ \vdots }{ \vdots }
      \TABROW{ 0 }{ 0 }{ \ldots }{ -\frac12w_1^{R,N-1}  }{ \frac12w_0^{R,N-1} + \frac32w_1^{R,N-1} }{ \frac12w_0^{R,N-1} }
      \TABROW{ 0 }{ 0 }{ \ldots }{ 0 }{ -\frac12w_1^{R,N} }{ \frac12w_0^{R,N} + \frac32w_1^{R,N} }
    \end{pmatrix}
    \nonumber
    \label{eqn::R}
  \end{align}
\end{small}
These definitions allows us to rewrite polynomial reconstruction as
the matrix-vector product:
\begin{align*}
  \hat{U}_L = \mathbf{L}\hat{U}
  \quad\textrm{and}\quad
  \hat{U}_R = \mathbf{R}\hat{U}.
\end{align*}

The reconstruction strategy can be easily extended to the
multidimensional case.
For example, neglect for simplicity of exposition the stochastic
dimensions and take $n=3$.
If $\mathbf{I}$ is the (properly sized) identity matrix and $\otimes$
denote the matrix Kronecker product, we have that
\begin{description}
\item[$(i)$] the reconstruction along the direction $x$ provides the left and
  right values at the cell interfaces parallel to the plane $y-z$ by applying the
  matrix operators
  $\mathbf{L}\otimes\mathbf{I}\otimes\mathbf{I}$ and
  $\mathbf{R}\otimes\mathbf{I}\otimes\mathbf{I}$;

\item[$(ii)$] the reconstruction along the direction $y$ provides the left and
  right values at the cell interfaces parallel to the plane $z-x$ by
  applying the matrix operators
  $\mathbf{I}\otimes\mathbf{L}\otimes\mathbf{I}$ and
  $\mathbf{I}\otimes\mathbf{R}\otimes\mathbf{I}$;
  
\item[$(iii)$] the reconstruction along the direction $z$ provides the left and
  right values at the cell interfaces parallel to $x-y$ by applying
  the matrix operators
  $\mathbf{I}\otimes\mathbf{I}\otimes\mathbf{L}$ and
  $\mathbf{I}\otimes\mathbf{I}\otimes\mathbf{R}$.
\end{description}
These formulas can be adapted in a straightforward way to build the reconstruction
matrices at the quadrature points, denoted by $\mathbf{Q}^{\beta}$, for $\beta \in \{ 1, \ldots, \hat{\mathcal{N}}\}$ with $\hat{\mathcal{N}}$ the number of quadrature points per dimension.  Then, along the $k$-th direction, irrespective of direction $k$ being
along a physical or a stochastic dimension, the application of the quadrature reconstruction matrices is given by:

\begin{align*}  
  \underbrace{  \mathbf{I}\otimes\ldots\otimes\mathbf{I}  }_{\text{$k-1$ times}}
  \otimes
  \underbrace{  \mathbf{Q}^\beta  }_{\text{$k$-th position}}
  \otimes
  \underbrace{\mathbf{I}\otimes\ldots\otimes\mathbf{I}}_{\text{$d-k$ times}}.
\end{align*}
The reconstruction matrices at the quadrature points are assembled in a nearly identical way to the $\mathbf{L}$ and $\mathbf{R}$ interface reconstruction matrices, but require differently scaled weights which depend upon the quadrature points.  See \cite{Shu} for more details.  The quadrature reconstruction matrices are necessary in both the physical and stochastic dimensions so that one may compute the cell averages and expectations as indicated by equation \eqref{eqn::general:update}.

However, a remarkable advantage of the tensor-train formulation of the SFV method, described in the next section, is that we can compute 1-dimensional reconstruction matrices and apply them core by core to the tensor-train representation without ever constructing a multi-dimensional analoque of \eqref{eqn::L} or \eqref{eqn::R}. This procedure is described in detail in algorithms \ref{alg::TT_SFV_PHYS_RECONSTRUCTION} and \ref{alg::TT_SFV_STOCH_RECONSTRUCTION}. 




\section{Tensor-train reformulation of the stochastic finite volume method}
\label{sec::TT-SFV-method}



We denote by $\ell=1,2,\ldots,p$ the indices for some sequence. The index $k$ runs over the range $1,2,\ldots,d$. Generalizing the discussion of the previous section, the variable $n$ takes on the values $1,2,3$, representing the number of directions of the physical space with $\bs{x} \in D_{\bs{x}}\subset \mathbb{R}^n$. The number of directions of the stochastic space is given by $m$, which is equal to the number of stochastic variables $\by\in D_{\by}$. The total number of dimensions is denoted by $d=n+m$ as before.

The scalar index $N$ is used if the physical space is one-dimensional, or alternatively, $N=(N_1,N_2,\ldots,N_n)$ represents a multi-index which counts the total number of partitions of the physical domain along every direction. Similarly, for the stochastic space, $M$ is a scalar index in one-dimensional cases, or $M=(M_1,M_2,\ldots,M_m)$ is a multi-index indicating the total number of partitions of the stochastic domain along every direction.

The index $i$ runs from $1$ to $N$ and labels the partitions of the one-dimensional physical space, or in multidimensional cases, $\idx{i}=(i_1,i_2,\ldots,i_n)^T$ with $i_s=1,2,\ldots,N_s$ for $s=1,2,\ldots,n$, serves as a multi-index labeling the partitions of the physical domain along each direction. The index $j$ labels the partitions of the one-dimensional stochastic space, or in the multidimensional context, $\idx{j}=(j_1,j_2,\ldots,j_m)^T$ where $j_s=1,2,\ldots,M_s$ for $s=1,2,\ldots,m$ is the corresponding multi-index.

We define $K_{\bs{x}}^i$ as a multidimensional cell in the physical space and $K_\by^j$ as a multidimensional cell in the stochastic space. The Cartesian product $\{K_{\bs{x}}^i\}\times\{K_\by^j\}$ constitutes a multidimensional mesh, as before.

The vector $\bv{U} = (U_1,U_2,\ldots,U_p)^T$ consists of $p$ components of conserved quantities $U_\ell(t,\bs{x},\by)$. The vector $\bv{\bar{U}}=(\bar{U}_{1},\bar{U}_{2},\ldots,\bar{U}_{p})^T$ comprises multidimensional arrays that collect the cell-averages over the mesh $\{K_{\bs{x}}^i\}\times\{K_\by^j\}$. Each array $U_\ell, \bar{U}_{\ell}$ for $\ell=1,2,\ldots,p$ is a multidimensional tensor.


\subsection{Tensor-train format}
\label{subsec::TT-intro} 
Let $d$ be a positive integer. By definition, a \emph{$d$-dimensional tensor}
${U}_\ell\in\mathbb{R}^{n_1\times n_2\times\ldots n_d}$ is an array of
real numbers whose entries are addressed by a $d$-sized multi-index
$i=(i_1,i_2,\ldots,i_d)$, e.g.,
${U}_\ell=\big({U}_\ell(\idx{i})\big)=\big({U}_\ell(i_1,i_2,\ldots,i_d)\big)$.
If $d=1$ or $d=2$, we prefer the usual terms \emph{vector} and
\emph{matrix}.
The integers $n_{k}$, $k=1,2,\ldots,d$, are called the
\emph{mode sizes} of the tensor.

We say that a $d$-dimensional tensor is in the \emph{tensor-train
format}~\cite{Oseledets-TT-2011}, denoted $\ten{U}_\ell$ if
there exist
two matrices ${\mathcal{U}}_{\ell,1}\in\mathbb{R}^{\nss{1}\times\rss{1}}$ and
${\mathcal{U}}_{\ell,d}\in\mathbf{R}^{\rss{d-1}\times\nss{d}}$
and $d-2$ three-dimensional tensors
${\mathcal{U}}_{\ell,k}\in\mathbb{R}^{r_{k-1}\times n_k \times r_k}$, $2\leq k\leq d-1$,
called the \emph{TT cores}, such that
\begin{align}
  \ten{U}_\ell(i_1,i_2,\ldots,i_d) =
  \sum_{\alpha_1=1}^{r_1}\sum_{\alpha_2=1}^{r_2}\ldots\sum_{\alpha_{d-1}=1}^{r_{d-1}}
  {\mathcal{U}}_{\ell,1}(i_1,\alpha_1){\mathcal{U}}_{\ell,2}(\alpha_1,i_2,\alpha_2)\ldots{\mathcal{U}}_{\ell,d}(\alpha_{d-1},i_d).
  \label{eqn:tt-format:def}
\end{align}
A tensor in tensor-train format will be referred to as a \emph{TT tensor} in the sequel.
The $(d-1)$ integer numbers $r_1,\ldots,r_{d-1}$ are the
\emph{TT ranks} of the TT tensor.
Introducing two further integers $\rss{0}=\rss{d}=1$, we can
restate~\eqref{eqn:tt-format:def} using only three-way tensors
$\ten{U}_{\ell,k}\in\mathbb{R}^{r_{k-1}\times n_k \times r_k}$, $1\leq
k\leq d$. Moreover, write  $\tenb{U} = (\ten{U}_{1},\ten{U}_{2},\ldots,\ten{U}_{p})^T$ to denote a $p$-sized vector of tensors of conserved quantities $\ten{U}_{\ell}(t)\in\mathbb{R}^{N_1\times N_2\times\ldots\times N_d}$ for $\ell=1,2,\ldots,p$ in tensor-train format.  

\begin{figure}[t!]
  \centering
  \includegraphics[width=0.65\textwidth]{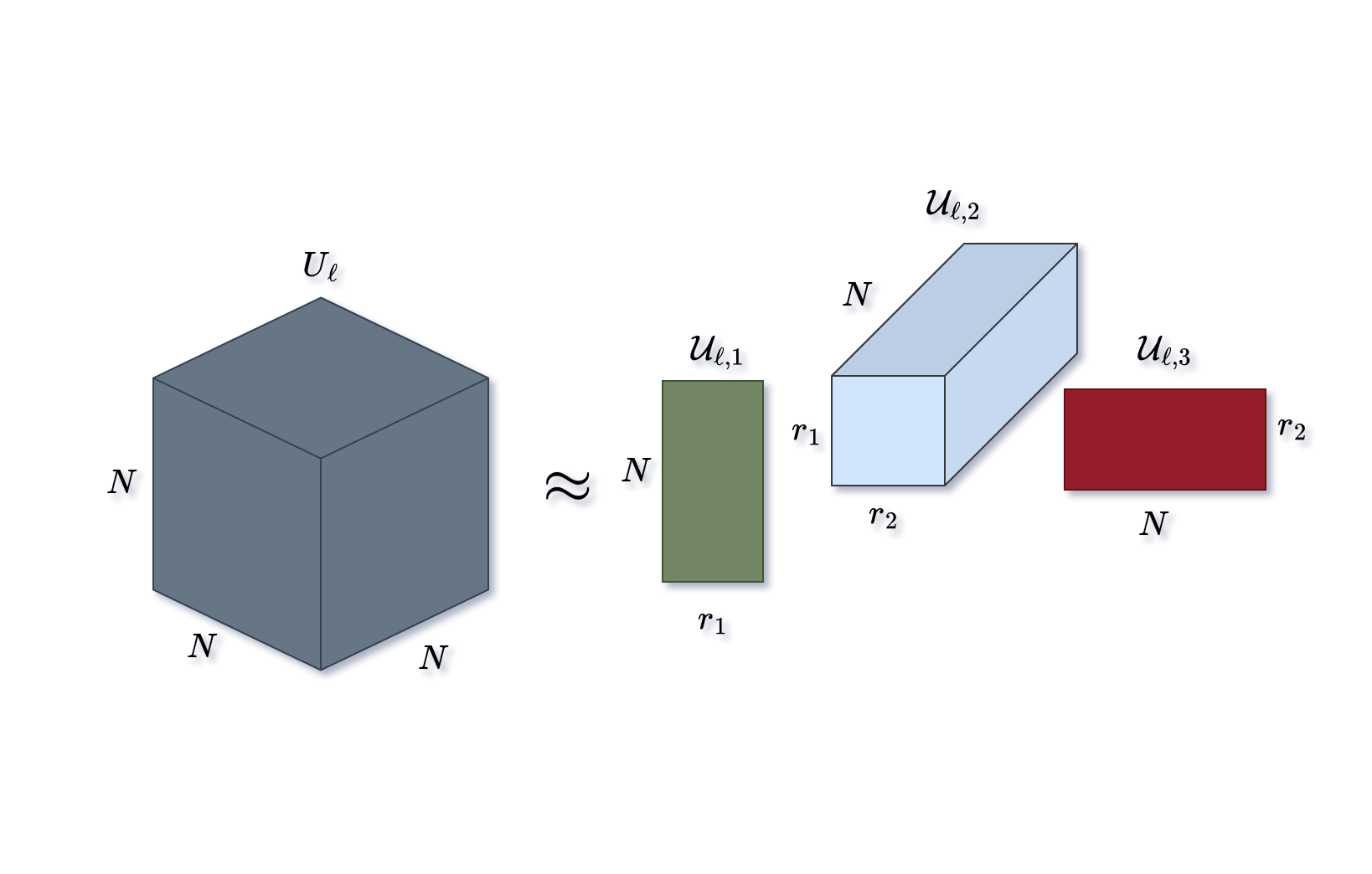}
  \caption{Graphical depiction of a three dimensional tensor $U_\ell$, representing a conserved quantity, and its tensor-train approximation via the tensor-train cores ${\mathcal{U}}_{\ell,k}$ for $k=1,2,3$,
    cf.~\eqref{eqn:tt-format:def}. Here,
    the rank set is $r=\big\{r_1,r_2\big\}$ and all cores have an equal number of modes
    $N=M$.}
  \label{fig:u_bar_TT_diagram}
\end{figure}

In Figure.~\ref{fig:u_bar_TT_diagram}, a three-dimensional tensor
$\bv{U}_{\ell}$ and its tensor-train representation, $\ten{U}_\ell$, with all mode sizes equal to $n$ and ranks
$r_1,r_2$ is depicted schematically via its cores, ${\ten{U}}_{\ell,k}$ for $k=1,2,3$.

The TT ranks $(\rss{1},\rss{2},\ldots,\rss{d-1})$ determine the size
of a TT tensor, i.e, the compression that the TT format
representation can achieve.
Let $\frak{n}=\max\limits_{k}\,\nss{k}$ and $\frak{r}=\max\limits_{k}\,\rss{k}$ be upper bounds of
the mode sizes $n_k$ and the ranks $r_k$.
The storage of a $d$-dimensional TT tensor $\ten{U}$ with mode sizes
$n_1,n_2,\ldots,n_d$ and ranks $r_0,r_1,r_2,\ldots,r_{d-1},r_d$ is
equal to $\sum_{k=1}^{d}n_kr_{k-1}r_{k}$, which is roughly 
proportional to $\mathcal{O}(d\frak{n}\frak{r}^2)$, see \cite{Oseledets-TT-2011}.
The linear dependence on $d$ makes the tensor-train representation
possibly free of the curse of dimensionality, and, if $\frak{r}\ll \frak{n}$, highly
compressive compared with the normal storage, which is proportional to
$\mathcal{O}(\frak{n}^d)$.

The TT format has several attractive properties.
First, consider the decomposition
\begin{align*}
  {U} = \ten{U} + \ten{E},
\end{align*}
where tensor $\ten{E}$ represent the approximation error.
Let $\hat{{U}}_k$ be the $k$-th matrix unfolding of ${U}$,
corresponding to the MATLAB/NumPy operation
\begin{align*}
  \hat{{U}}_k = \text{reshape}\Big( {U}, \big[ \Pi_{\ell=1}^k\nss{k}, \Pi_{\ell=k+1}^k\nss{d} \big]\Big),
\end{align*}
and set $\rss{k}=\RANK{(\hat{{U}}_k)}$ for $k=1,2,\ldots,d-1$.
An exact tensor-train decomposition exists with these ranks, so that
$\norm{\ten{E}}{\mathcal{F}}=0$, cf.~\cite[Theorem~2.1]{Oseledets-TT-2011}.
If we set $\rss{k}\leq\RANK(\hat{{U}}_k)$ for $k=1,2,\ldots,d-1$, then,
the low-rank \emph{best approximation} in the Frobenius norm with
these ranks, e.g., ${U}_{\text{best}}$, exists,
cf~\cite[Corollary~2.4]{Oseledets-TT-2011}, and the algorithm TT-SVD
\cite{Oseledets-TT-2011} provides a \emph{quasi-optimal} approximation to
it, in the sense that
\begin{align*}
  \norm{ {U} - \ten{U} }{\mathcal{F}}
  \leq
  \sqrt{d-1}\norm{ {U} - {U}_{\text{best}} }{\mathcal{F}}.
\end{align*}
If the truncation tolerance for the singular value decomposition (SVD)
of each unfolding is set to $\delta = \frac{\epsilon}{\sqrt{d-1}}
\norm{\ten{U}}{\mathcal{F}}$, the algorithm TT-SVD is able to
construct a tensor-train approximation $\ten{{U}}$ such that 
\begin{align*}\label{eqn::tt_err}
  \norm{ U - \ten{U} }{\mathcal{F}} \leq \epsilon \norm{ {U} }{\mathcal{F}}.
\end{align*}

The TT-format allows us to store tensors in a compact form and to
perform fundamental arithmetic operations efficiently without leaving
the format, cf.~\cite{Oseledets-TT-2011}.
For example, the element-wise (Hadamard) product $\circ$ of two
tensors $\ten{U}$ and $\ten{V}\in\mathbb{R}^{n_1\times
  n_2\times\dots\times n_d}$ with TT-ranks $r_k\leq \frak{r}=\max_{k}(r_k)$
and $s_k\leq \frak{s}=\max_{k}(s_k)$, respectively, is a tensor
$\ten{W}=\ten{U}\circ\ten{V}$ with TT-ranks at most $r_k s_k$
and can be computed in $\mathcal{O}(d\frak{n}\frak{r}^2\frak{s}^2)$ operations,
see~\cite{Oseledets-TT-2011}.
To compute a TT-decomposition of a given tensor we can apply a
generalization of the matrix Singular Value Decomposition (SVD) to the
case of multi-dimensional arrays presented in
Ref.~\cite{Oseledets-TT-2011}, which we will refer to as the
\emph{TT-SVD algorithm}.
This algorithm is numerically stable and has a complexity proportional
to $\mathcal{O}( d\frak{n}^2\frak{r}^3 )$.
Moreover, according to~\cite[Theorem~2]{Oseledets-TT-2011}, an exact
TT representation always exists, although not unique, and an
approximate TT representation with a prescribed accuracy and minimal
rank can be found.
Other tensor-train decomposition algorithms such as those based on the
\emph{cross interpolation} makes it possible to construct an
approximate representation of a tensor using only a subset of its
entries, see Refs.\cite{OSELEDETS201070,savostyanov2011fast}.
Low-rank cross interpolation relies on the matrix skeleton
decomposition~\cite{Goreinov1997}.
The skeleton decomposition of a matrix $A \in \mathbb{R}^{p \times q}$
is defined as
\begin{align*}
  \mat{A}' = \mat{A}(:, \mathcal{J}) \mat{A} (\mathcal{I}, \mathcal{J})^{-1} \mat{A}( \mathcal{I}, :),
\end{align*}
where $\mathcal{I} = (i_1, \ldots, i_r)$ and $\mathcal{J} = (j_1,
\ldots, j_r)$ are subsets of the index sets $\big[1,\ldots, p\big]$ and
$\big[1,\ldots, q\big]$.
The optimal submatrix $\mat{A}(\mathcal{I},\mathcal{J})$ is determined
by maximizing the determinant modulus among all $r\times r$
submatrices of $\mat{A}$.
Finding such a matrix is a computationally challenging
task~\cite{BARTHOLDI1982190}.
However, an approximation can be obtained using the \emph{maxvol
algorithm}~\cite{max_vol}, which is computationally efficient and requires
$2c(n - r)r$ operations, where $c$ is typically a small constant in
practical applications.
The selection of indices $(\mathcal{I},\mathcal{J})$ carried by this
algorithm is intended to capture most of the information in $\mat{A}$
through the decomposition.

\subsection{Tensor-train reformulation of the SFV method}
\label{subsec::SFV_TT}
Our next task is to reformulate the SFV method
\eqref{eqn::general:update} in the TT format.
Let ${\bar{U}}_{\ell}(t)$, $\ell=1,2,\ldots,p$,denote the multidimensional array collecting the cell averages of the
$\ell$-th component of $\bv{U}(t,\bs{x},\bv{y})$ at some given instant
$t\geq0$ on the multidimensional grid with dimensions $d=n+m$.
We denote its tensor-train representation as $\ten{\bar{U}}_{\ell}$,
and we number its TT cores by adding a subindex $k$,
i.e., $\ten{\bar{U}}_{\ell,k}$, $k=1,2,\ldots,d$.

We construct the cell-average tensor of the initial solution by
cell-wise integration of the cores of $\ten{U}_{\ell}$ as described
in Algorithm \ref{alg::TT_SFV_CELL_INT}, which is inspired by \cite{DOLGOV2020106869, Alexandrov-Manzini-Skau-Truong-Vuchkov:2023} and \cite{OSELEDETS201070} and adapted to our purposes.
This algorithm requires sampling the initial solution
$U_{\ell}(0,\bs{x},\by)$ at a set of cell quadrature nodes, thus resulting in
a TT tensor $\ten{U}_{\ell}(0)$ that can be directly
built in the tensor-train format by using the cross approximation
algorithm.
The separation of variables due to the TT representation allows us to
carry out the numerical integration fiber by fiber and independently across cores, thus
providing a very efficient implementation, to obtain the TT tensor represention of \eqref{eqn::init_Ubar}, $\bar{\ten{U}}(0)$.
\begin{algorithm}[!htb]
  \caption{TT-SFV Algorithm: cell-wise integration to compute $\ten{\bar{U}}_{\ell}$, $\ell=1,2,\ldots,p$
    \label{alg::TT_SFV_CELL_INT}}
  \begin{algorithmic}[1]
    \Let{$\ten{{U}}_{\ell}$}{$U_{\ell}(0,x,y)$} \Comment{Compute $\big(\ten{{U}}_{\ell}\big)$ using cross interpolation sampling $U_{\ell}(0,x,y)$ at quadrature points}
    \State{Perform cell-wise integration along fibers to obtain a core list $(\ten{\bar{U}}_{\ell,k})$ for the cell average state}
    \For{$k = 1 \textrm{ to } d$}
    \For{$j = 1: \mathcal{\hat{N}}:\mathcal{\hat{N}}$}
    \Let{$\ten{\bar{U}}_{\ell,k}(:,j/\mathcal{\hat{N}},:)$}{$\sum_{s=0}^{\mathcal{\hat{N}}} \ten{U}_{\ell,k}\omega(j+s)$}
    \EndFor
    \EndFor
    \State{Assemble $\ten{\bar{U}}_{\ell}$ from the core list as in equation \eqref{eqn:tt-format:def}. }
  \end{algorithmic}
\end{algorithm}

From the cell averaged state $\ten{\bar{U}}_{\ell}$ we reconstruct
the approximate solution trace at the cell interfaces using the WENO
reconstruction techniques.
This operation is performed by multiplying each physical core
$\ten{\bar{U}}_{\ell,k}$, $k=1,2,\ldots,n$ by the left and right
global reconstruction matrices along the appropriate physical dimension.
\begin{algorithm}[!h]
  \caption{TT-SFV Algorithm: Physical Reconstruction 
    \label{alg::TT_SFV_PHYS_RECONSTRUCTION}}
  \begin{algorithmic}[1]
    \Let{$\ten{\bar{U}}_{\ell,s}$}{$\ten{\bar{U}}_{\ell}$}\Comment{For the desired physical dimension, $s$, obtain the corresponding TT core.}
    \State{Compute $\mathbf{L}$ and $\mathbf{R}$ as given in equations \eqref{eqn::L} and \eqref{eqn::R}.}
    \State{Compute $\mathbf{L}\ten{\bar{U}}_{\ell,s}$ and $\mathbf{R}\ten{\bar{U}}_{\ell,s}$.}
    \For{$j = 0 \textrm{ to } n$ and $j\neq{s}$}
    \If{$s=1$}
    \Let{$\ten{\bar{U}}_{\ell,L}^1$}{$\mathbf{L}\ten{\bar{U}}_{\ell,1}\,\Big(\prod_{k=2}^{d} \ten{\bar{U}}_{\ell,k}\Big)$}
    \Let{$\ten{\bar{U}}_{\ell,R}^1$}{$\mathbf{R}\ten{\bar{U}}_{\ell,1}\,\Big(\prod_{k=2}^{d} \ten{\bar{U}}_{\ell,k}\Big)$}
    \Else
    \Let{$\ten{\bar{U}}_{\ell,L}^s$}{$\Big(\prod_{k=1}^{s-1} \ten{\bar{U}}_{\ell,k}\Big)\,\mathbf{L}\ten{\bar{U}}_{\ell,s}\,\Big(\prod_{k=i+1}^{d} \ten{\bar{U}}_{\ell,k}\Big)$}
    \Let{$\ten{\bar{U}}_{\ell,R}^s$}{$\Big(\prod_{k=1}^{s-1} \ten{\bar{U}}_{\ell,k}\Big)\,\mathbf{R}\ten{\bar{U}}_{\ell,s}\,\Big(\prod_{k=i+1}^{d} \ten{\bar{U}}_{\ell,k}\Big)$}
    \EndIf
    \EndFor
  \end{algorithmic}
\end{algorithm}
Applying such matrices computes the reconstructed values at the
quadrature points in every dimension, and numerical integration fiber by fiber
returns the tensor-train representation of the numerical fluxes.
\begin{figure}[!hb]
  \centering
  \includegraphics[width=0.65\textwidth]{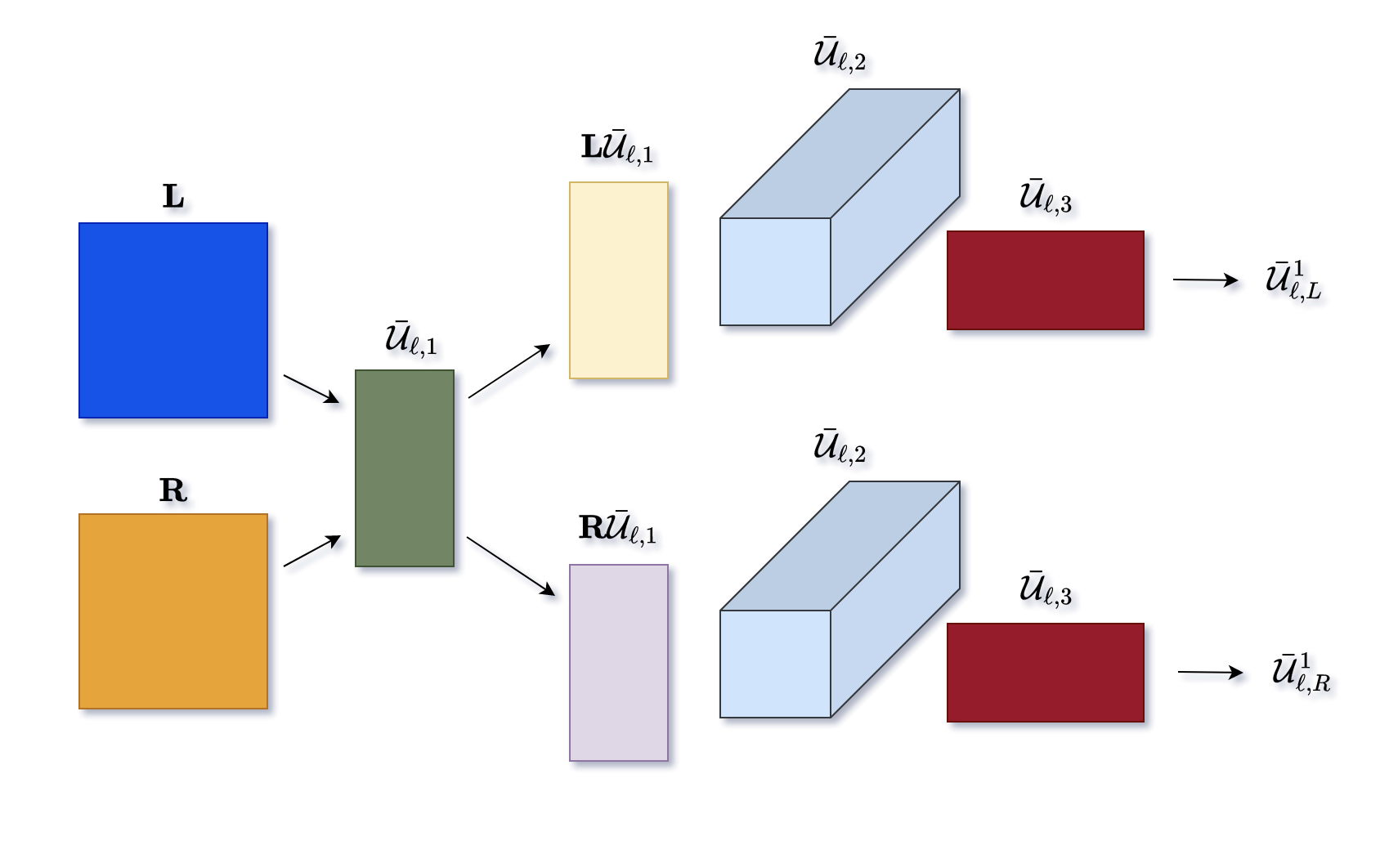}
  \caption{Graphical representation of applying reconstruction in
    physical domain for the case $n=1$ and $m=2$.
    The left and right reconstruction matrices are applied to the
    core of $\bar{\mathcal{U}}_{\ell}$ representing the physical dimension, $\bar{\mathcal{U}}_{\ell,1}$, resulting in reconstruction along
    the modes.
    The result is $ \mathbf{L}\bar{\mathcal{U}}_{\ell,1}$ and $ \mathbf{R}\bar{\mathcal{U}}_{\ell,1}$ which
    are then recombined with the stochastic cores to produce
   $\bar{\mathcal{U}}^1_{\ell,L}$ and $\bar{\mathcal{U}}^1_{\ell,R}$, as desired. 
  }
    \label{fig:physical_reconstruction_TT_diagram}
\end{figure}
This procedure for the reconstruction in the physical space is
outlined in Algorithm \ref{alg::TT_SFV_PHYS_RECONSTRUCTION} and
depicted schematically in Diagram
\ref{fig:physical_reconstruction_TT_diagram} for the special case
$n=1$ and $m=2$, giving $3$ total dimensions.  The resulting TT tensors, with reconstructions in physical direction $s$, is denoted by $\ten{U}_{\ell,L}^s$ and $\ten{U}_{\ell,R}^s$ for the left and right reconstructions, respectively.

\begin{figure}[!ht]
  \centering
  \includegraphics[width=0.65\textwidth]{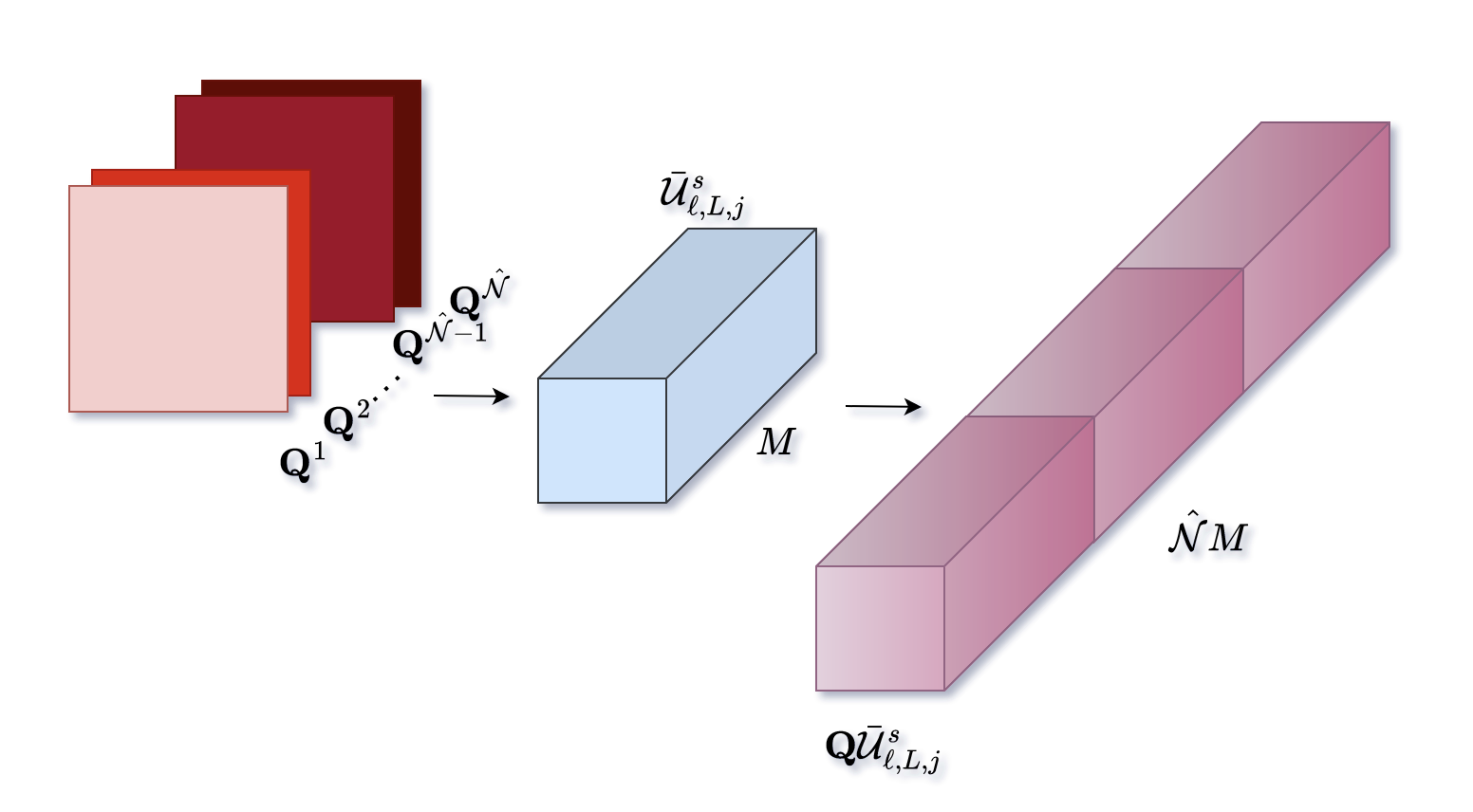}
  \caption{Graphical representation of performing reconstruction to the
    stochastic core, $j$, of the left reconstructed tensor of quantity $\ell$ in physical direction $s$ with $\hat{\mathcal{N}}$ quadrature points.
    The matrices $Q^\beta$ for $\beta = 1,2,\ldots, \hat{\mathcal{N}}$ are the reconstruction
    matrices for quadrature points.
    Each quadrature reconstruction matrix is applied to each slice of the stochastic core from, e.g., $ \bar{\mathcal{U}}^s_{\ell,L,j}$ and stored cell-wise in $\mathbf{Q}\bar{\mathcal{U}}^s_{\ell,L,j}$,
    which scales the number of stochastic modes, $M$, by the number of quadrature points.  If quadrature reconstruction at a physical core $j\neq s$ with mode length $N$ is desired, the procedure is identical.
    }
  \label{fig:stochastic_reconstruction_TT_diagram}
\end{figure}

The procedure to obtain the reconstruction at the quadrature points
for the stochastic dimensions is detailed in
Algorithm \ref{alg::TT_SFV_STOCH_RECONSTRUCTION}.  The application of quadrature reconstruction matrices to a particular stochastic core is illustrated in Diagram
\ref{fig:stochastic_reconstruction_TT_diagram}.  As can be seen from Algorithm \ref{alg::TT_SFV_STOCH_RECONSTRUCTION}, we make copies of the left and right reconstructed tensors $\bar{\mathcal{U}}_{\ell, L}^s$, $\bar{\mathcal{U}}_{\ell, R}^s$ for a specific physical direction $s$. 
We then proceed recursively by running over the cores corresponding to stochastic dimensions, $k = n +1 , \ldots d$.  In the case $n > 1$, an identical reconstruction procedure must be performed on the cores corresponding to the physical dimension $j\in\{1,\ldots, n\}$, $j\neq s$.  The resulting TT tensors are denoted by $\ten{U}^s_{\ell, RecL}$ and $\ten{U}^s_{\ell, RecR}$ each with resulting mode sizes $N=(\hat{\mathcal{N}}N_1, \hat{\mathcal{N}}N_2,\ldots,N_s, \ldots, \hat{\mathcal{N}}N_n)$ and $M=(\hat{\mathcal{N}}M_1, \hat{\mathcal{N}}M_2,\ldots, \hat{\mathcal{N}}M_m)$.
\begin{algorithm}[!ht]
  \caption{TT-SFV Algorithm: Stochastic Reconstruction\label{alg::TT_SFV_STOCH_RECONSTRUCTION}}
  \begin{algorithmic}[1]
    \State{ Assume there is at least one physical dimension, $n \geq 1$. Let $s$ denote the physical dimension under consideration.}
    \Let{$\ten{\bar{U}}_{L_{qp}}$}{$\ten{\bar{U}}_{\ell, L}^s$}\Comment{Copy the left and right reconstructed TT tensors of physical dimension $s$.}
    \Let{$\ten{\bar{U}}_{R_{qp}}$}{$\ten{\bar{U}}^s_{\ell, R}$}
    \For{$k = n+1 \textrm{~to~} d $}
    \For{ $\alpha_{k-1} = 1 \textrm{~to~} r_{k-1}$ }
    \For{ $\alpha_{k}  = 1 \textrm{~to~} r_{k}$ }
    \For{$\beta = 1 \textrm{ to } \hat{\mathcal{N}}$}
    \State{Construct ${\mathbf{Q}^\beta}$ for $\ten{\bar{U}}_{L_{qp},k}(\alpha_{k-1},:,\alpha_{k})$
     and $\ten{\bar{U}}_{R_{qp},k}(\alpha_{k-1},:,\alpha_{k})$.}
    \State{ $\mathbf{C}^\beta_L = \mathbf{Q}^\beta\ten{\bar{U}}_{L_{qp},k}(\alpha_{k-1},:,\alpha_{k})$}\Comment{Temporary storage vectors.}
    \State{ $\mathbf{C}^\beta_R = \mathbf{Q}^\beta\ten{\bar{U}}_{R_{qp},k}(\alpha_{k-1},:,\alpha_{k})$}
    \For{$\gamma = 1 \textrm{ to } \mathcal{M}$}
    \State{ $\mathbf{Q}\ten{\bar{U}}_{L_{qp},k}(\alpha_{k-1}, \gamma(\hat{\mathcal{N}}) + \beta ,\alpha_k) = \mathbf{C}^\beta_L(\gamma) $ }\Comment{ Resulting mode length is $\hat{\mathcal{N}}\mathcal{M}$.}
    \State{ $\mathbf{Q}\ten{\bar{U}}_{R_{qp},k}(\alpha_{k-1}, \gamma(\hat{\mathcal{N})} + \beta ,\alpha_k) = \mathbf{C}^\beta_R(\gamma) $}
    \EndFor
    \EndFor
    \EndFor
    \EndFor
    \State{Re-build $\ten{\bar{U}}_{L_{qp}}$ and $\ten{\bar{U}}_{R_{qp}}$ from quadrature reconstructed cores.}
    \If{$k=d$}
    \Let{$\ten{\bar{U}}_{L_{qp}}$}{$\Big(\prod_{j=1}^{d-1} \ten{\bar{U}}_{L_{qp},j}\Big)\,\mathbf{Q}\ten{\bar{U}}_{L_{qp},k}$}
    \Let{$\ten{\bar{U}}_{R_{qp}}$}{$\Big(\prod_{j=1}^{d-1} \ten{\bar{U}}_{R_{qp},j}\Big)\,\mathbf{Q}\ten{\bar{U}}_{R_{qp},k}$}
    \Else
    \Let{$\ten{\bar{U}}_{L_{qp}}$}{$\Big(\prod_{j=1}^{k-1} \ten{\bar{U}}_{L_{qp},j}\Big)\,\mathbf{Q}\ten{\bar{U}}_{L_{qp},k}\,\Big(\prod_{j=k+1}^{d} \ten{\bar{U}}_{L_{qp},j}\Big)$}
    \Let{$\ten{\bar{U}}_{R_{qp}}$}{$\Big(\prod_{j=1}^{k-1} \ten{\bar{U}}_{R_{qp},j}\Big)\,\mathbf{Q}\ten{\bar{U}}_{R_{qp},k}\,\Big(\prod_{j=k+1}^{d} \ten{\bar{U}}_{R_{qp},j}\Big)$}
    \EndIf
    \EndFor
    \State{After final recursion, assign to $\bar{\mathcal{U}}_{\ell, RecL}^s$ and $\bar{\mathcal{U}}_{\ell, RecR}^s$ }
    \Let{$\bar{\mathcal{U}}_{\ell, RecL}^s$}{$\ten{\bar{U}}_{L_{qp}}$}
    \Let{$\bar{\mathcal{U}}_{\ell, RecR}^s$}{$\ten{\bar{U}}_{R_{qp}}$}
    
  \end{algorithmic}
\end{algorithm}

With the reconstruction procedure complete, we evaluate the fluxes using the cross-approximation algorithm.  Let $\hat{\bv{F}}^{TT}\big(\cdot, \cdot\big)$ denote the function which uses the cross-approximation algorithm (see the preceding subsection) to approximate the element-wise evaluation of the numerical flux of equation \eqref{eq:laxfriedrichsflux}. We will refer to $\hat{\bv{F}}^{TT}\big(\cdot, \cdot\big)$ as the TT flux.  For a specific physical dimension $s$, the arguments passed in to the TT flux are the left and right fully reconstructed tensors $\bar{\mathcal{U}}_{\ell, RecL}^s,\bar{\mathcal{U}}_{\ell, RecR}^s$. 
To obtain the TT approximation of $\bar{\bv{F}}_{i\pm\half,j}(t)$ of equation \eqref{eq:numflux:def}, we introduce the shift operators $\mathcal{T}^\pm_s$.  The action of $\mathcal{T}^\pm_s$ applied to a TT tensor is simply to shift the entries of the fibers in core $s$ to the next index ($+$) or the preceding index ($-$).  Then, for physical direction $s$ we obtain $\hat{\bv{F}}^{TT}_{s+\frac12}\big(\bar{\mathcal{U}}_{\ell, RecL}^s,\mathcal{T}^+_s\bar{\mathcal{U}}_{\ell, RecR}^s\big)$ and $\hat{\bv{F}}^{TT}_{s-\frac12}\big(\mathcal{T}^+_s\bar{\mathcal{U}}_{\ell, RecL}^s,\bar{\mathcal{U}}_{\ell, RecR}\big)$.  Note carefully that the subscript $s$ indicates a physical direction in contrast to $\bar{\bv{F}}_{i\pm\half,j}(t)$ which has the subscript $i$ to denote a cell index.  The notation $s\pm\frac12$ can then be understood as accounting for all cell indices $i$ in the direction $s$.  With these definitions in hand, it is easy to define the derivative approximation in direction $s$, for which we write 
\begin{equation}
    \delta_s \hat{\bs{F}}^{TT} = \frac{1}{\Delta x_s}\Big(\hat{\bv{F}}^{TT}_{s+\frac12}\big(\bar{\mathcal{U}}_{\ell, RecL}^s,\mathcal{T}^-_s\bar{\mathcal{U}}_{\ell, RecR}^s\big) -  \hat{\bv{F}}^{TT}_{s-\frac12}\big(\mathcal{T}^+_s\bar{\mathcal{U}}_{\ell, RecL}^s,\bar{\mathcal{U}}_{\ell, RecR}\Big),
\end{equation}
with $\Delta x_s$ the vector of cell lengths in the direction $s$. The cell-wise averages in the physical directions $j\in\{1,\ldots, n\}$, $j\neq s$ are computed using Algorithm \ref{alg::TT_SFV_CELL_INT}.  The expectations are also computed using Algorithm \ref{alg::TT_SFV_CELL_INT} but with the additional step of applying the measure $\mu(\by)$.  To do this a TT representation of $\mu(\by)$, denoted by $\mu(\by)^{TT}$. It is easy to obatin the stochastic cell measures $| K_\by^j|$ from the core-wise integration algorithm applied to the cores of $\mu(\by)^{TT}$.  Then, to compute the expectation of the $j$th core, we simply perform an element wise product of the $j$th core of $\delta_s \hat{\bv{F}}^{TT}$ and $\frac{\mu(\by_j)^{TT}}{|K_\by^j|}$, written $\hat{\bv{F}}^{TT}_j \odot \frac{\mu(\by_j)^{TT}}{|K_\by^j|}$, and the cell-wise integration Algorithm \ref{alg::TT_SFV_CELL_INT} is applied as before.  The result of this procedure, detailed in Algorithm \ref{alg::TT_SFV_Flux}, is a list of cores with which we can compute $\delta_s\bar{\bs{F}}^{TT}$.

We may now state the fully discrete TT-SFV method. Again uniformly partitioning the time integration interval as $0=t^{0}<\ldots<t^{k}<t^{k+1}<\ldots<t^{K}=T$ and using the forward Euler update for concreteness, we may write the fully discrete update for a quantity $\ell$ as
\begin{equation}
    \bar{\ten{U}}_\ell^{k+1} = \bar{\ten{U}}_\ell^{k} - \Delta t^k \sum_{s=1}^{n}\delta_s\bar{\bs{F}}^{TT,k}.
\end{equation}
Here, for simplicity, the source terms are assumed to be zero. In the case of nonzero source terms, one can construct them exactly as $\bar{\ten{U}}$.  As previously mentioned, we use a second- or third-order SSP time integration method.  After each stage of the time integration, rounding is applied to control the size of the ranks of the TT approximation \cite{Oseledets-TT-2011}.
\begin{algorithm}[!htb]
  \caption{TT-SFV Algorithm: Numerical Flux
    \label{alg::TT_SFV_Flux}}
  \begin{algorithmic}[1]
    \Let{$\hat{\bv{F}}^{TT}_{s-\frac12}$}{$\hat{\bv{F}}\Big(\mathcal T^-_s\bar{\mathcal{U}}_{\ell, RecL}^s,\bar{\mathcal{U}}_{\ell, RecR}^s\Big)$}\Comment{ For a desired physical dimension $s$, apply left shift operator, $\mathcal T^{-}_s$, to the right reconstructed state and send with the left reconstructed state as arguments to numerical flux via cross approximation to obtain $\hat{\bv{F}}^{TT}_{s-\frac12 }$.}
    \Let{$\hat{\bv{F}}^{TT}_{s+\frac12}$}{$\hat{\bv{F}}\Big(\bar{\mathcal{U}}_{\ell, RecL}^s, \mathcal T^+_s\bar{\mathcal{U}}_{\ell, RecR}^s\Big)$}\Comment{ For a desired physical dimension $s$, apply right shift operator, $\mathcal T^{+}$, to the left reconstructed state and send with the right reconstructed state as arguments to numerical flux via cross approximation to obtain $\hat{\bv{F}}^{TT}_{i+\frac12 }$.}
    \State{$\delta_s \hat{\bv{F}}^{TT} = \frac{1}{\Delta x_s}\Big(\hat{\bv{F}}^{TT}_{s+\frac12 }  - \hat{\bv{F}}^{TT}_{s-\frac12 }\Big)$}\Comment{Compute derivative in $i$th physical dimension}
    \For{$j = 1 \textrm{ to }  n$ and $j \neq s$ }
    \Let{$\delta_s\hat{\bv{F}}_j^{Int}$}{$\int \delta_s \hat{\bv{F}}^{TT}_j \mathrm d x_j$} \Comment{Use Algorithm \ref{alg::TT_SFV_CELL_INT} to perform cell-wise integration in $j$th core of $\delta_s \hat{\bv{F}}^{TT}$. }
    \EndFor
    \For{$j = n+1 \textrm{ to }  d$ }\Comment{Use Algorithm \ref{alg::TT_SFV_CELL_INT} to compute $j$th expectation. }
    \Let{$\mathbb E_j\Big[\delta_i\hat{\bv{\bv{F}}}^{TT}\Big]$}{$\frac{1}{|K_\by^j|}\int \delta_s \hat{\bv{F}}^{TT}_j \odot \mathrm d \mu(\by_j)^{TT}$} 
    \EndFor
    \If{$s=1$}\Comment{Assemble $\delta_s \bar{\bv{F}}^{TT}$ from the cores computed above. }
    \Let{$\delta_s \bar{\bv{F}}^{TT}$}{$\delta_1 \hat{\bv{F}}^{TT}_{1}\Big(\prod_{k=2}^n \delta_1\hat{\bv{F}}_k^{Int}\Big) \Big(\prod_{\ell = n+1}^{d} \mathbb E_\ell\Big[\delta_1\hat{\bv{F}}^{TT}\Big]\Big) $}
    \Else
    \Let{$\delta_s \bar{\bv{F}}^{TT}$}{$\Big(\prod_{k=1}^{s-1}\delta_s\hat{\bv{F}}_k^{Int}\Big)\delta_s \hat{\bv{F}}^{TT}_{s}\Big(\prod_{k=s+1}^n \delta_s\hat{\bv{F}}_k^{Int}\Big) \Big(\prod_{\ell = n+1}^{d} \mathbb E_\ell\Big[\delta_s\hat{\bv{F}}^{TT}\Big]\Big) $}
    \EndIf
  \end{algorithmic}
\end{algorithm}

In the next section, we provide several numerical tests for the TT-SFV method with which to validate the procedures discussed in the foregoing discussion.  Two of the tests can be found elsewhere in the literature \cite{AbgrallTokareva-2017, sfv-adaptive-SIAM, HertyKolbMuller2023}, which are low-dimensional with $n=m=1$.  These serve as a sanity check for the method.  The remaining tests show the application of the TT-SFV method on problems with $n\in\{1,2,3\}$ and $m$ varying between $8$ and $12$.



\section{Numerical results}
\label{sec::tests}
We present several test problems with which to verify the methods
presented above.
In each case we specify a domain $ I\times
D_{\bs{x}}\times D_{\by}\subset\mathbb{R}^+\times\mathbb{R}^n\times\mathbb{R}^m$,
where $I$ is the time interval, $D_{\bs{x}}$ is the spatial domain and $D_{\by}$
is the stochastic domain.

For our experiments, we use the tntorch library
\cite{usvyatsov2022tntorch} which is based on the pytorch machine
learning framework \cite{pytorch} or the TT-toolbox \cite{GitHub:Oseledets-TT-Toolbox}, though other TT libraries/toolboxes
\cite{torchtt, GitHub:Oseledets-TTPY}
could have just as easily been used.
To use the TT-SFV method, one only needs to be able to represent a tensor in TT
format, use the cross-approximation algorithm and apply the rounding
operation.  
We did not exploit any GPU acceleration or parallelization and all
tests were run on a personal Macbook Pro.


\subsection{Stochastic Linear Advection: Smooth Initial Condition}
\begin{figure}[!h]
  \centering
  \includegraphics[width=0.6\textwidth]{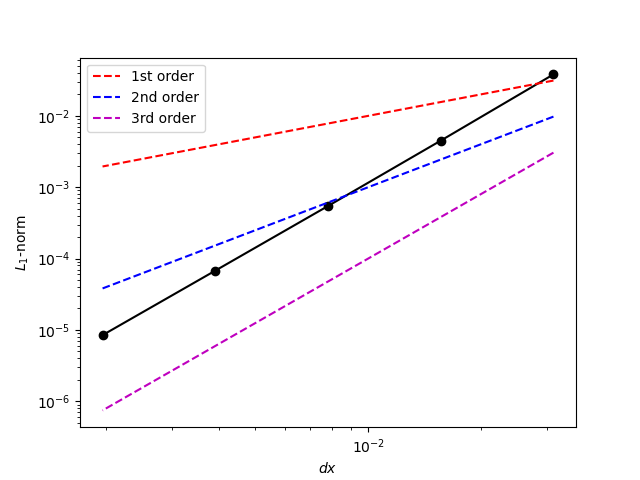}
  \caption{Results from a h-refinement convergence study for the problem \eqref{eqn::smooth_linear_advection} demonstrating that the SFV in TT format obtains the expected 3rd order convergence.}
  \label{fig::smooth_lin_transport}
\end{figure}

Here we test the convergence properties of the SFV method in TT format on a smooth linear advectin problem with $m=n=1$ so that  $I = (0,0.1]$, $D_x=[0,1]$ and $D_y = [0,1]$ and the one dimensional stochastic variables are assumed to be uniformly distributed. Thus, we have 
\begin{equation}\label{eqn::smooth_linear_advection}
  \begin{aligned}
    \partial_t u(t,x,y)
    + a \, \partial_{x} u(t,x,y)
    = 0 && (t,x,y) \in I\times D_x \times D_y \subset  \mathbb R^+ \times \mathbb R \times \mathbb R,\\
    u(0,x,y) = u_0(x,y) && (t,x,y) \in \{ 0 \} \times D_x \times D_y,
  \end{aligned}
\end{equation}
with the initial condition, following \cite{AbgrallTokareva-2017},
\begin{equation}\label{eqn::smooth_linear_advection_ic}
  u_0(x,y) = \sin(2\pi(x+0.1y)).
\end{equation}
As demonstrated in \cite{AbgrallTokareva-2017}, the SFV method with WENO3 reconstruction in the physical and stochastic space should obtain 3rd order convergence in the expectation as the grid is refined.  We verify in figure \ref{fig::smooth_lin_transport} that the optimal convergence rates are indeed obtained.  For this test, we set the TT approximation threshold, see equation \eqref{eqn::tt_err}, to $\epsilon = 10^{-6}$, and employ the SSP(3,3) time integration scheme \cite{Gottlieb2009} with rounding (subsection \ref{subsec::TT-intro}) applied after each stage. 
The initial grid has $2^5$ grid points and is succesively refined by factors of $2$ until we reach $2^{10}$. For every refinement level, the ranks are never greater than $4$ as the simulation evolves.  

Here, we are merely interested in verifying that the SFV method in TT-format is capable of reproducing previous results by obtaining theoretically predicted convergence rates.  Given the low-dimensionality of the problem, it is reasonable to expect that the TT-format will provide comparable runtimes to a traditional SFV code.

\subsection{Stochastic Inviscid Burgers'}
Here we solve the inviscid Burgers' equation with $n=m=1$, giving  
\begin{equation}
  \begin{aligned}
    \partial_t u(t,x,y)
    + \frac{1}{2}\partial_x u^2(t,x,y) = 0 && (t,x,y) \in I\times D_x \times D_y\subset\mathbb{R}^+\times\mathbb{R}\times\mathbb{R},\\
    u(0,x,y) = u_0(x,y) && (t,x,y) \in \{ 0 \} \times D_x \times D_y,
  \end{aligned}
  \label{eqn::Burgers}
\end{equation}
  where the initial condition is smooth, defined by 
\begin{equation}
  u_0(x,y) = \sin(2\pi x)\sin(2\pi y),
  \label{eqn::Burgers_ic}
\end{equation}
where we take  $D_x=[0,1]$ and $D_y = [0,1]$ with $y\sim \mathrm{Unif}(D_y)$ or $y\sim \mathrm{Beta}(2,5)$.  The same problem can be found in 
\cite{HertyKolbMuller2023} and \cite{sfv-adaptive-SIAM}.  As in those works, we take $I = [0, 0.35]$ for our time interval.

\begin{figure}[!h]
  \centering
  \includegraphics[width=0.8\textwidth]{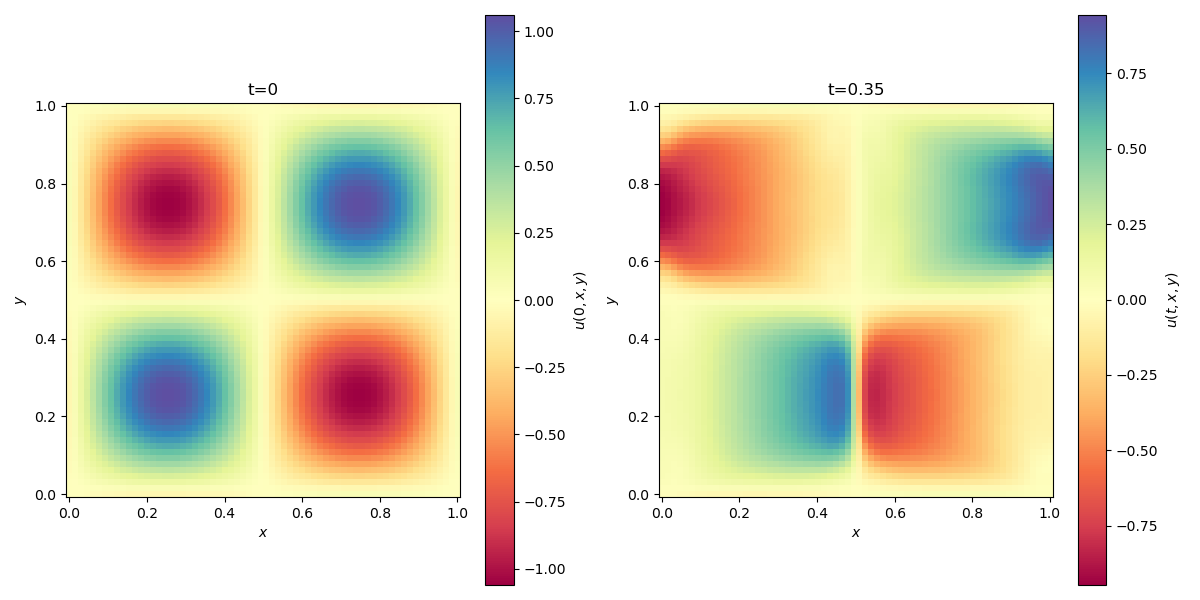}
  \caption{Stochastic Burgers'. Solution to the problem \eqref{eqn::Burgers}. The initial condition is shown on the left and the TT-SFV solution at t=0.35 is shown on the right.}
  \label{fig:stoch_Burgers}
\end{figure}

The solution to the problem is shown in figure \ref{fig:stoch_Burgers} with the initial condition for comparison. The stationary shocks which develop due to the nonlinearity of the problem are clearly captured in the physical and stochastic directions.

\begin{figure}[!h]
  \centering
  \includegraphics[width=0.48\textwidth]{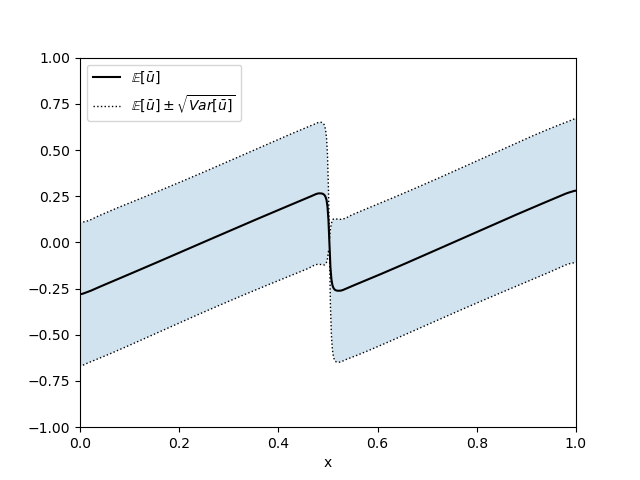}
  \includegraphics[width=0.48\textwidth]{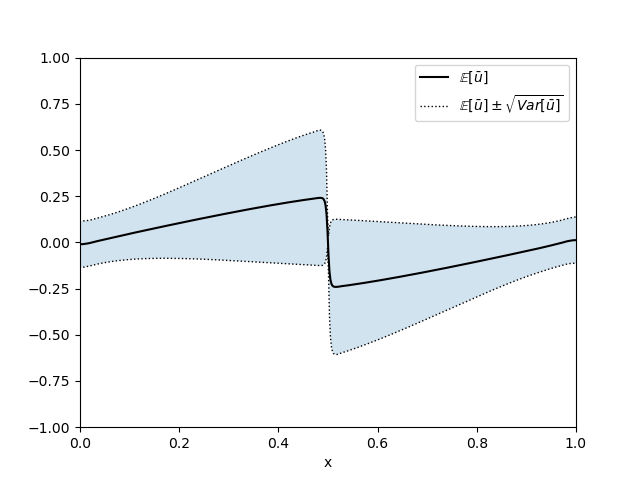}
  \caption{Stochastic Burgers'. Expectation values with standard deviations shown the problem \eqref{eqn::Burgers} with uniformly distributed random variables on the left and Beta distributed random variables on the right.  The solutions are in good agreement with the previous works.}
  \label{fig:stoch_Burgers_E}
\end{figure}

The expectation values at the final time $t=0.35$ are given in figure \ref{fig:stoch_Burgers_E} for both distributions.  It is clear that the TT-SFV algorithm is capable of handling the stationary shock as well as picking up the skewed Beta(2,5) distribution shown on the right.  The figure shows that our results are in good agreement with the results found in \cite{HertyKolbMuller2023, sfv-adaptive-SIAM}.  
For the results shown, we set $\epsilon = 9e-04$, below the error of the WENO3 SFV discretization, and set the grid points to $512$. The maximal rank observed as the simulation progressed was $10$ but usually $4$, in both cases.  As in the previous case, given the low-dimensionality of the problem, we did not compare runtimes against a traditional SFV code as there is no expected gain for the TT-format in this low-dimensional setting.  The solution was computed with the SSP(3,3) time integration scheme.

\subsection{A Stochastic Sod Problem}
TT-SFV method is now applied to the Euler equations for a Sod-like test
problem.  Here, we take $n=1$ and $m=12$ with $I = (0,0.2]$, $D_x=[0,1]$, and $D_\by = [0,1]^{12}$.
The problem is constructed so that the initial condition takes
stochastic inputs as in the previous examples.
\begin{align}
  & \partial_t \bv{U}(t,x,\by)+\partial_x \bv{F}\big(\bv{U}(t,x,\by)\big)=\bv{0}, \quad (t,x,\by) \in I\times D_x \times D_\by \subset  \mathbb R^+ \times \mathbb R \times \mathbb R^{12}, \label{eq:1D-PDE:stochastic} \\
  & \bv{U}(0,x,\by) = \bv{U}_0(x,\by) \quad (0,x,\by) \in \{0\} \times D_x \times D_\by.  
\end{align}
with $\by \sim \textrm{Unif}\big(D_\by\big)$ and
\begin{equation}
  \bv{U}(t,x,\by) = \left[\begin{array}{c} \rho \\ \rho u \\ \rho e \end{array}\right](t,x,\by), \qquad 
  \bv{F}\big(\bv{U}(t,x,\by)\big) = \left[\begin{array}{c} \rho u\\ \rho u^2 + p \\ u (\rho e + p) \end{array}\right](t,x,\by),
\end{equation}
where the pressure is defined by 
\begin{equation}
  p(t,x,\by) = (\gamma-1)(\rho e-\frac12\rho u^2).
\end{equation}
The parameterized random initial data is given by
\begin{equation}
  \begin{aligned}
    \left[\begin{array}{c} \rho_0 \\ \quad \\  \rho_0 u_0 \\ \quad \\  p_0 \end{array}\right] =
    \begin{cases}
      1.0 + 0.1\by_1 - 0.05\by_7, &\quad \text{if} \quad x < 0.5, \\
      0.125-0.05\by_2+0.1\by_8,   &\quad \text{if} \quad x > 0.5. \\
      0.05\by_3-0.01\by_9,        &\quad \text{if} \quad x < 0.5 \\
      0.05\by_4-0.01\by_{10},      &\quad \text{if} \quad x > 0.5. \\
      1.+0.1\by_5-0.05\by_{11},    &\quad \text{if} \quad x < 0.5, \\
      0.1+0.05\by_6-0.01\by_{12},  &\quad 
      \text{if} \quad x > 0.5
    \end{cases}
  \end{aligned}
\end{equation}
In addition, we assume that $\gamma = 1.4$, though it could also be
made to depend upon an additional random parameter.

\begin{figure}[!hb]
  \centering
  \includegraphics[width=0.65\textwidth]{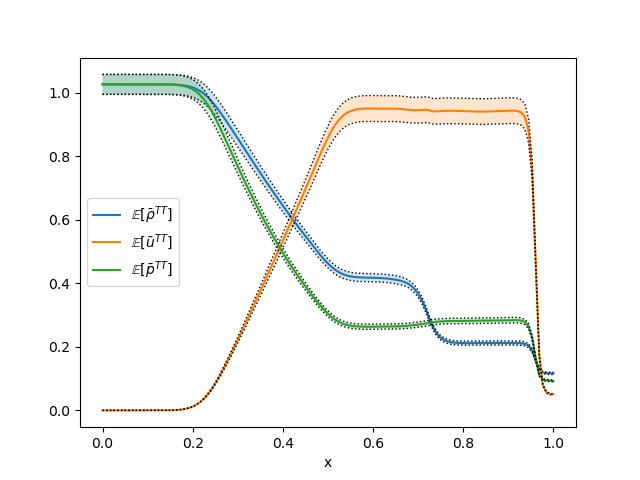}
  \caption{Test~Case~3: Computed expectations for $\bar{\rho}^{TT}$,
    $\bar{u}^{TT}$ and $\bar{p}^{TT}$ for the Sod problem with one
    physical and twelve stochastic dimension and $128$ cells per
    dimension.
    Here $t=0.2$.
    Black dotted lines indicate plus or minus the standard deviation for the respective variable, as in previous plots.
  }
  \label{fig:Sod_1x_12s}
\end{figure}

In figure \ref{fig:Sod_1x_12s}, the expectations with standard deviations are shown for the density, velocity and pressure at $t=0.2$.  The solution was computed with $128$ cells in each of the twelve dimensions and the SSP(2,2) time integration scheme.  For the depicted solution we set the ranks to $r=r_k=1$ for $k\in\{0,\ldots,d\}$, providing a significant speed up which we investigate further in the next subsection.

\subsection{Stochastic Shock Bubble}
\begin{figure}[!h]
  \centering
  \includegraphics[width=0.45\textwidth]{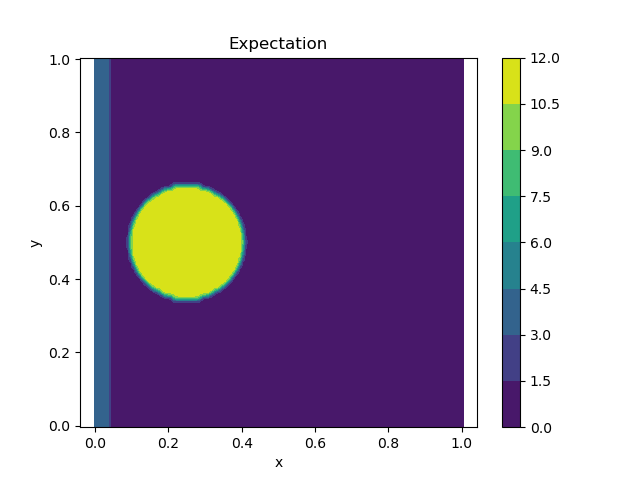}
  \includegraphics[width=0.45\textwidth]{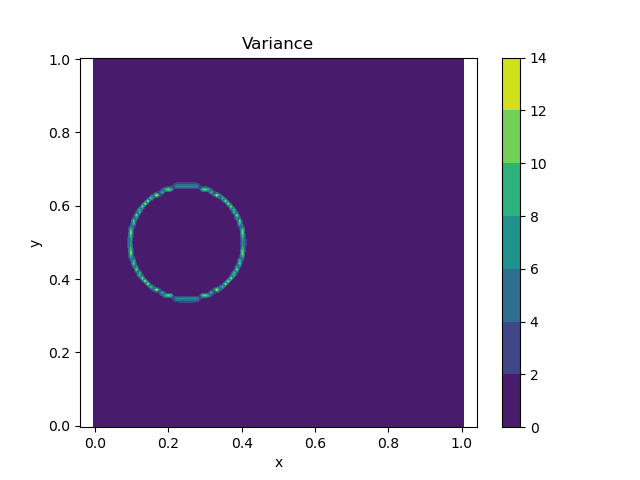}\\
  \includegraphics[width=0.45\textwidth]{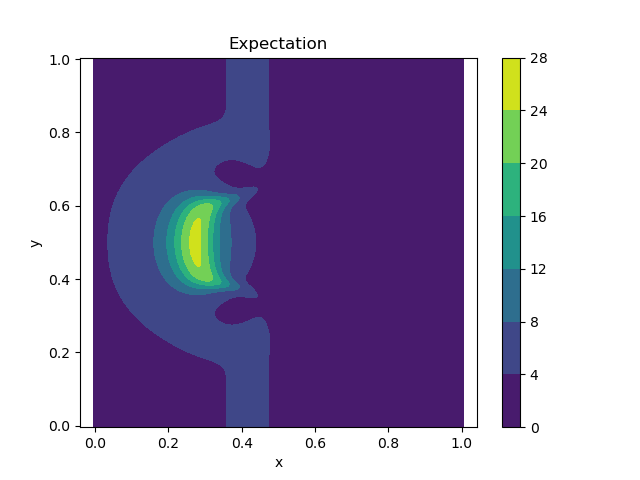}
  \includegraphics[width=0.45\textwidth]{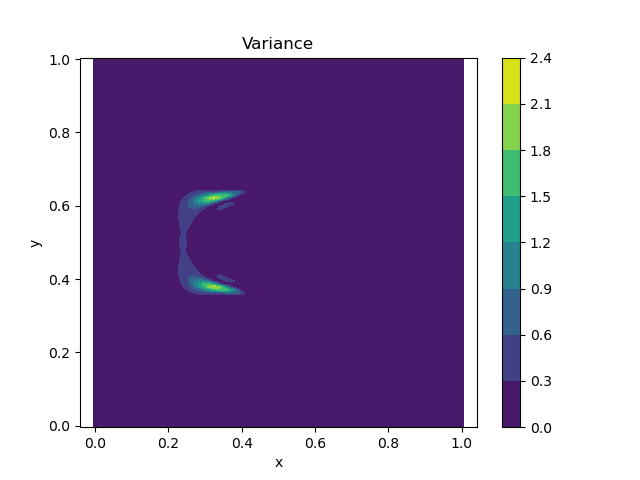}
  \caption{Stochastic Shock Bubble.  The expectation (left) and variance (right) of the solution to \eqref{eqn:shock_bubble} with $N=M=128$.  The top row depicts the initial conditions and the bottom row shows the evolution at time $t=0.03$.
  }
  \label{fig:shock_bubble_2x_12s}
\end{figure}

We again solve the Euler equations, now with $n=2$ and $m=12$ and $I = (0,0.2]$, $D_{\bs{x}}=[0,1]^2$, $D_\by = [0,1]^{12}$ and $\by \sim U\big(D_\by\big)$ . The problem reads
\begin{equation}\label{eqn:shock_bubble}
\begin{aligned}
  & \partial_t \bv{U}(t,\bs{x},\by)+ \nabla_{\bs{x}} \cdot \bv{F}\big(\bv{U}(t,\bs{x},\by)\big)=\bv{0}, \quad (t,\bs{x},\by) \in I\times D_{\bs{x}} \times D_\by \subset  \mathbb R^+ \times \mathbb R^2 \times \mathbb R^{12},  \\
  & \bv{U}(0,\bs{x},\by) = \bv{U}_0(\bs{x},\by) \quad (0,\bs{x},\by) \in \{0\} \times D_{\bs{x}} \times D_\by.
\end{aligned}
\end{equation}
with
\begin{equation}
  \bv{U}(t,\bs{x},\by) = \left[\begin{array}{c} \rho \\ \rho \bv{u} \\ \rho e \end{array}\right](t,\bs{x},\by), \qquad 
  \bv{F}\big(\bv{U}(t,\bs{x},\by)\big) = \left[\begin{array}{c} \rho \bv{u}\\ \rho \bv{u}^2 + p \\ \bv{u} (\rho e + p) \end{array}\right](t,\bs{x},\by),
\end{equation}
where $\bv{u} = ( u, v)^T$. The initial data is given by
\begin{equation}
  \begin{aligned}
    \left[\begin{array}{c} \rho_0 \\ \quad \\  \rho_0 u_0 \\ \quad \\ \rho_0 v_0 \\ \quad \\  p_0 \end{array}\right] =
    \begin{cases}
      3.86859+0.1\by_4-0.05\by_{11} + 0.5B(\bs{x},\by)  &\quad \text{if} \quad \bs{x}_1 < 0.04, \\
      1 +0.05\by_5-0.01\by_{12} +  0.5B(\bs{x},\by),   &\quad \text{if} \quad \bs{x}_1 > 0.04, \\
      11.2536+\by_7,        &\quad \text{if} \quad \bs{x}_1 < 0.04, \\
      0.01\by_8,      &\quad \text{if} \quad \bs{x}_1 > 0.04, \\
      0,        &\quad \text{if} \quad \bs{x}_1 < 0.04, \\
      0,      &\quad \text{if} \quad \bs{x}_1 > 0.04, \\
      167.345+10\by_9,    &\quad \text{if} \quad \bs{x}_1 < 0.04, \\
      1+0.1\by_{10},  &\quad 
      \text{if} \quad \bs{x}_2 > 0.04,
    \end{cases}
  \end{aligned}
\end{equation}
where $B(\bs{x},\by)$ is defined by
\begin{equation}
  \begin{aligned}
    B(\bs{x},\by) =
    \begin{cases}
      10 + 0.1\by_6  &\quad \text{if} \quad \sqrt{(\bs{x}_1-0.25)^2+(\bs{x}_2-0.5)^2} < 0.15 + 0.01\by_1, \\
      0,   &\quad \text{if} \quad \sqrt{(\bs{x}_1-0.25)^2+(\bs{x}_2-0.5)^2} > 0.15 + 0.01\by_1,
    \end{cases}
  \end{aligned}
\end{equation}
where the function $B(\bs{x},\by)$ introduces uncertainty in the bubble radius. The density function is quite challenging for the cross-approximation algorithm since it is not a simple low-rank function.  To this end, the density approximation is allowed to have a higher rank in the first few cores than the other variables. We set the ranks for the TT approximation of the density to $r = [32,32,32,1,1,1,1,1,1,1,1,1,1,1]$. Thus, even with higher rank cores in the approximation, the attained compression is quite satisfactory.  We show the expectation and variance of the solution in figure \ref{fig:shock_bubble_2x_12s}.  The solution was computed with $N=128$ and the SSP(2,2) time integration scheme.  We intend to investigate further the coupling of full rank solution with low-rank stochastic representations in a future work.

\subsection{Scaling Study and Speed Up}

\begin{figure}[!ht]
  \centering
  \includegraphics[width=0.65\textwidth]{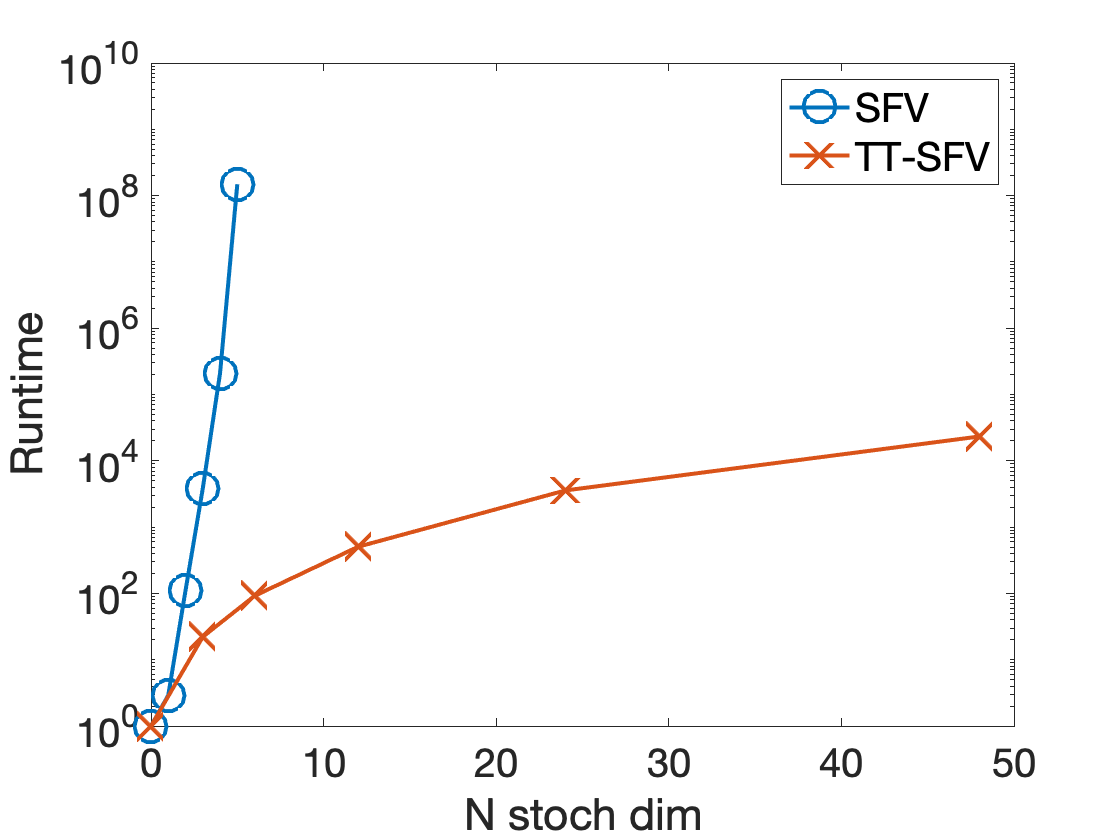}
  \caption{Test~Case~3: Runtimes in seconds for Sod problem with
    increasing number of stochastic dimensions (purple crosses and
    yellow line) against the C++ implementation (red circles and blue
    line).
  }
  \label{fig:Sod_scaling}
\end{figure}

To test runtime scaling with dimension, we run a Sod problem with the
physical dimension fixed at one and increase the stochastic dimension
for successive runs.
Further, the runtimes of the nonlinear TT-SFV algorithm are compared
with a C++ implementation of the SFV
method~\cite{AbgrallTokareva-2017,SFVSpringer2014}. The results of the Sod test are shown in Fig.~\ref{fig:Sod_scaling}.  We see that the runtimes are comparable until the stochastic dimension is greater than three.  Around $5$ stochastic dimension the runtimes for the C++ code begin to become intractable.  We remark that for this test both codes are run in serial on a CPU.





\section{Conclusion}\label{sec::conclusion}
In this paper, we have addressed the challenge of applying the
stochastic finite volume (SFV) method to hyperbolic systems with
uncertainties, particularly those involving shock waves and
discontinuities. Our integration of the SFV method with the
tensor-train (TT) decomposition has resulted in a novel TT-SFV method
that significantly mitigates the curse of dimensionality. The global
weighted essentially non-oscillatory (WENO) reconstruction within the
TT framework allows for accurate representation of discontinuous
solutions without introducing spurious oscillations.

Through our numerical experiments, including the Sod shock tube
problem under various stochastic dimensions, we have demonstrated the
efficiency, accuracy, and stability of the TT-SFV method. The results
show a promising reduction in computational complexity while
maintaining a high level of accuracy, even as the number of stochastic
dimensions increases. Moreover, the comparative analysis between our
Python implementation and traditional C++ SFV code indicates that our
method is not only theoretically sound but also practically viable.

Despite the advancements presented in this work, several limitations
and avenues for future research remain. The current TT-SFV method may
require further refinement for handling more complex boundary
conditions and for extending to higher dimensions in physical space. Future
work will also focus on improving the efficiency of the algorithm,
optimizing the TT-rank selection process, and exploring adaptive
strategies to dynamically adjust the TT-rank during simulations based
on the evolving solution features.

In conclusion, the TT-SFV method represents a substantial step forward
in the field of uncertainty quantification for hyperbolic systems. It
opens the door to solving a broader class of problems where
traditional methods are limited by computational constraints. This
work not only provides a new tool for researchers but also sets the
stage for future developments in the numerical treatment of stochastic
PDEs with discontinuities.




\section*{Acknowledgements}
The Laboratory Directed Research and Development (LDRD) program of
Los Alamos National Laboratory partially supported the work of
GM under project number 2023ZZZZDR, and ST and SW under project number
20220121ER.
Los Alamos National Laboratory is
operated by Triad National Security, LLC, for the National Nuclear
Security Administration of U.S. Department of Energy (Contract
No.\ 89233218CNA000001).
G. Manzini is a member of the Italian Gruppo Nazionale Calcolo
Scientifico - Istituto Nazionale di Alta Matematica (GNCS-INdAM).
The manuscript is approved for public release as Technical Report
LA-UR-24-21630


\bibliographystyle{siamplain}
\bibliography{references}

\end{document}